\documentclass[review]{article}
\usepackage[latin9]{inputenc}
\usepackage{units}
\usepackage{amsmath}
\usepackage{amssymb}
\usepackage{graphicx}
\usepackage{esint}
\usepackage[numbers]{natbib}

\makeatletter
\@ifundefined{showcaptionsetup}{}{%
 \PassOptionsToPackage{caption=false}{subfig}}
\usepackage{subfig}
\makeatother

\begin{document}

\title{Investigation of a Structured Fisher's Equation with Applications
in Biochemistry}

\author{John T. Nardini\thanks{Department of Applied Mathematics, University of Colorado, Boulder
80309-0526, United States (john.nardini@colorado.edu, dmbortz@colorado.edu)}, D.M. Bortz$^{*}$}
\maketitle
\begin{abstract}
Recent biological research has sought to understand how biochemical
signaling pathways, such as the mitogen-activated protein kinase (MAPK)
family, influence the migration of a population of cells during wound
healing. Fisher's Equation has been used extensively to model experimental
wound healing assays due to its simple nature and known traveling
wave solutions. This partial differential equation with independent
variables of time and space cannot account for the effects of biochemical
activity on wound healing, however. To this end, we derive a structured
Fisher's Equation with independent variables of time, space, and biochemical
pathway activity level and prove the existence of a self-similar traveling
wave solution to this equation. We also consider a more complicated
model with different phenotypes based on MAPK activation and numerically
investigate how various temporal patterns of biochemical activity
can lead to increased and decreased rates of population migration.
\end{abstract}
\textbf{keywords:} Stage-structure, Traveling Wave Solutions, Wound Healing, Biochemical Signaling Pathways

\section{Introduction\label{sec:Introduction}}

Traveling wave solutions to partial differential equations (PDEs)
are often used to study the collective migration of a population of
cells during wound healing \citep{cai_multi-scale_2007,denman_analysis_2006,landman_diffusive_2005,landman_travelling_2007,maini_traveling_2004},
tumorigenesis \citep{kuang_data-motivated_2015}, and angiogenesis
\citep{pettet_model_1996,sherratt_new_2001}. R.A. Fisher introduced
what is now referred to as Fisher's Equation in 1937 to model the
advance of an advantageous gene in a population \citep{fisher_wave_1937}.
Since then, it has been used extensively in math biology literature
to model the migration of a monolayer of cells during experimental
wound healing assays \citep{cai_multi-scale_2007,jin_reproducibility_2016,maini_traveling_2004}. 

Fisher's Equation is written as
\begin{equation}
u_{t}=Du_{xx}+\lambda u(K-u)\label{eq:fisher}
\end{equation}
with subscripts denoting differentiation with respect to that variable
and $u=u(t,x)$ representing a population of cells over time $t$
at spatial location $x$. The first term on the right hand side of
\eqref{eq:fisher} represents diffusion in space with rate of diffusion,
$D$, and the second term represents logistic growth of the population
with proliferation rate, $\lambda,$ and carrying capacity, $K$.
As shown in \citep[$\S$ 11.2]{murray_mathematical_2002}, \eqref{eq:fisher}
admits traveling wave solutions of the form
\[
u(t,x)=U(z),\ z=x-ct
\]
where $c$ denotes the speed of the traveling wave solution and $U(z)$
denotes the traveling wave profile. Traveling wave solutions to \eqref{eq:fisher}
thus maintain a constant profile, $U(z),$ over time that moves leftward
if $c<0$ or rightward if $c>0$ with speed $|c|$. It is also shown
that \eqref{eq:fisher} has a positive and monotonic profile for $|c|\ge2\sqrt{D\lambda}$,
which is biologically relevant when $u(t,x)$ denotes a population
of cells. Kolmogoroff proved in 1937 that any solution to \eqref{eq:fisher}
with a compactly-supported initial condition will converge to a traveling
wave solution with minimum wavespeed $c=2\sqrt{D\lambda}$ \citep{kolmogoroff_etude_1937}.
See \citep[$\S$ 5.4]{murray_lectures_1977} for a proof of this. There
is also a wide literature on studies into extensions of Fisher's Equation,
such as Fisher's Equation coupled with chemotaxis \citep{ai_reaction_2015,landman_diffusive_2005},
time-dependent rates of proliferation and diffusion \citep{hammond_analytical_2011},
and space-dependent rates of diffusion \citep{curtis_propagation_2012}.

Structured population models, or PDE models with independent variables
to distinguish individuals by some continuously-varying properties,
were first investigated via age-structured models in the early 20th
century \citep{mckendrick_applications_1927,sharpe_problem_1911}.
The 1970s saw a revival in structured population modeling after the
introduction of methods to investigate nonlinear structured population
models \citep{gurtin_non-linear_1974}, which led to our current understanding
of semigroup theory for linear and nonlinear operators on Banach spaces
\citep{webb_population_2008}. Several recent biological studies have
demonstrated the existence of traveling wave solutions to structured
population models \citep{ducrot_travelling_2011,ducrot_travelling_2009,gourley_vector_2007,so_reaction-diffusion_2001},
and another study used an independent variable representing subcellular
$\beta$-catenin concentration to investigate how signaling mutations
can cause intestinal crypts to invade healthy neighboring crypts \citep{murray_modelling_2010}.

Recent biological research has focused on the influence of biochemical
signaling pathways on the migration of a population of cells during
wound healing. Particular emphasis has been placed on the mitogen-activated
protein kinase (MAPK) signaling cascade, which elicits interesting
patterns of activation and migration in response to different types
of cytokines and growth factors in various cell lines \citep{chapnick_leader_2014,matsubayashi_erk_2004}.
For example, experimental wounding assays of madine darby canine kidney
cells (MDCKs) in \citep{matsubayashi_erk_2004} yielded a transient
pulse of ERK 1/2 (a specific MAPK protein) activity in the cell sheet
that only lasted for a few minutes. This pulse of activity was followed
by a slow wave of activity that propagated from the wound margin to
submarginal cells over the course of several hours. The second wave
was determined to be crucial for regulating MDCK sheet migration.
The authors of \citep{matsubayashi_erk_2004} proposed that these
fast and slow waves of ERK 1/2 activity could be caused by the production
of reactive oxygen species (ROS) and epidermal growth factor (EGF),
respectively. Similar experiments with fibroblasts also demonstrated
this first transient wave of ERK 1/2 activity, but not the following
slow wave. The authors of \citep{chapnick_leader_2014} found that
human keratinocyte (HaCaT) cells exhibit ERK 1/2 activity primarily
at the wound margin during similar experimental wound healing assays
with a high density in response to treatment with transforming growth
factor-$\beta$ (TGF-$\beta$).

In this study, we detail an approach to investigate a structured version
of Fisher's Equation that is motivated by the above experimental observations.
Previous structured population models have been restricted to traits
that primarily increase over time, such as age or size, but our analysis
allows for both activation and deactivation along the biochemical
activity dimension. 

In Section \ref{sec:Model-Development}, we develop our structured
population model and devote Section \ref{sec:Characteristic-equations-from}
to a review of relevant material from size-structured population models.
We demonstrate the existence of self-similar traveling wave solutions
to the model in Section \ref{sec:Traveling-waves-analysis}. We then
study a more complicated version of our model where migration and
proliferation of the population depend on MAPK activity levels in
Section \ref{sec:Numerical-simulations} before making final conclusions
and discussing future work in Section \ref{sec:Discussion-and-future}.

\section{Model Development\label{sec:Model-Development}}

We model a cell population during migration into a wound, denoted
by $u(t,x,m)$, for 
\[
u:[0,\infty)\times\mathbb{R}\times\left[m_{0},m_{1}\right]\rightarrow\mathbb{R}
\]
where $t$ denotes time, $x$ denotes spatial location, and $m$ denotes
activation along a biochemical signaling pathway with minimum and
maximum levels $m_{0}$ and $m_{1}$, respectively. As a first pass,
we assume that any cells of the same MAPK activity level will activate
identically over time in the same environment. This assumption allows
us to model the activation distribution of the population over time
deterministically by considering how cells of all possible MAPK activity
levels activate and deactivate over time. We note that biochemical
signaling is an inherently heterogeneous process, so our approach
would benefit from a further investigation with stochastic differential
equations. 

As discussed in \citep{de_roos_gentle_1996}, crucial aspects of a
structured population model include the individual state, the environmental
state, external forcing factors, and feedback functions. The \emph{individual
state} is a dimension used to distinguish between individuals of a
population and is typically based on physiological properties such
as age or size. As activation of biochemical signaling pathways influences
cell migration through diffusive and proliferative properties of cells,
we incorporate the biochemical activity dimension, $m$, as an individual
state for our model. 

The \emph{environmental state} of a population is the external factors
that influence individual behavior. Recall that external cytokines
and growth factors, such as ROS, TGF-$\beta$, and EGF, influence
activation of the MAPK signaling cascade and promote migration during
wound healing. The cell population will not directly affect the level
of external growth factor in this work, so an \emph{external forcing
factor} will be used to represent treatment with these chemicals here.
The external chemical concentration at time $t$ will be denoted by
$s(t),$ and the activation response of cells to this chemical will
be given by the function $f(s).$

A \emph{feedback function }included in our work will be the inhibition
of individual cell proliferation in response to a confluent density.
As proliferation is hindered by contact inhibition, we introduce a
new variable, 
\begin{equation}
w(t,x):=\int_{m_{0}}^{m_{1}}u(t,x,m)dm\label{eq:w_def}
\end{equation}
to represent the population of cells at time $t$ and spatial location
$x$. Proliferation of the population will accordingly vanish as $w(t,x)$
approaches the carrying capacity, $K$.

Our model, which we term as a \emph{structured Fisher's Equation},
is given by the PDE:
\begin{eqnarray}
u_{t}+\underbrace{\left(f(s(t))g(m)u\right)_{m}}_{activation} & = & \underbrace{D(m)u_{xx}}_{dif\negthinspace fusion}+\underbrace{\lambda(m)u\left(K-w(t,x)\right)}_{population\ growth}\label{eq:structured_fisher_eqn}\\
w(t,x) & = & \int_{m_{0}}^{m_{1}}u(t,x,m)dm\nonumber \\
u(t=0,x,m) & = & \phi(x,m)\nonumber \\
u(t,x,m=m_{1}) & = & 0\nonumber \\
w(t,-\infty)=K &  & w(t,x=+\infty)=0\nonumber 
\end{eqnarray}
 The function $g(m)\in C^{1}\left([m_{0},m_{1}]\right)$ denotes the
rate of biochemical activation in the population, $s(t)\in L^{\infty}(\mathbb{R}^{+})$
denotes the external chemical concentration in the population, $f(s)\in L_{loc}^{1}(0,\infty)$
denotes the activation response of cells to the level of signaling
factor present, $D(m)$ and $\lambda(m)$ denote biochemically-dependent
rates of cell diffusion and proliferation, and $\phi(x,m)$ denotes
the initial condition of $u$. The spatial boundary conditions specify
that the cell density has a confluent density at $x=-\infty$ and
an empty wound space at $x=+\infty.$ We use a no flux boundary condition
at $m=m_{1}$ so that cells cannot pass this boundary. In the remainder
of this study, we will write $f(s(t))$ as $f(t)$ for simplicity,
though we note that this function will differ between cell lines that
respond differently to the same chemical during wound healing\footnote{Note that an extension for modeling the dynamics governing $s(t)$
will be considered in a future study.}.

The solution space of \eqref{eq:structured_fisher_eqn}, $\mathcal{D}$,
is defined with inspiration from \citep{webb_population_2008} and
\citep[$\S$ 1.1]{volpert_traveling_1994}. If we let $Z$ denote the
space of bounded and twice continuously differentiable functions on
$\mathbb{R}$, then we define
\[
\mathcal{D}:=\left\{ u(t,x,m)\left|\int_{m_{0}}^{m_{1}}u(t,x,m)dm\in Z\right.\right\} ,
\]
i.e., $u(t,x,m)\in\mathcal{D}$ if $w(t,x)\in Z$ for all $t>0.$
We note that $\int_{m_{0}}^{m_{1}}\phi(x,m)dm$ need only be bounded
and piecewise continuous with a finite number of discontinuities \citep{volpert_traveling_1994}.
If $\phi(x,m)$ is not sufficiently smooth in $m,$ we obtain generalized
solutions of \eqref{eq:structured_fisher_eqn} \citep{webb_population_2008}. 

In Section \ref{sec:Traveling-waves-analysis}, we will investigate
\eqref{eq:structured_fisher_eqn} with constant rates of diffusion
and proliferation (i.e., $D(m)=D,\lambda(m)=\lambda$) and $f(t)=1$.
By substituting
\begin{eqnarray*}
u^{*}=u/K,\ \  &  & t^{*}=\lambda Kt,\\
x^{*}=x\sqrt{\lambda K/D},\ \  &  & m^{*}=(m-m_{0})/(m_{1}-m_{0}),\\
g^{*}(m^{*})= &  & \ g(m^{*}(m_{1}-m_{0})+m_{0})/(\lambda Km_{1}),
\end{eqnarray*}
and dropping asterisks for simplicity, \eqref{eq:structured_fisher_eqn}
can be non-dimensionalized to
\begin{eqnarray}
u_{t}+\underbrace{(f(t)g(m)u)_{m}}_{activation} & = & \underbrace{u_{xx}}_{diffusion}+\underbrace{u\left(1-\int_{0}^{1}u(t,x,m)dm\right)}_{population\ growth}\label{eq:structured_fisher_eqn_norm}\\
w(t,x) & = & \int_{0}^{1}u(t,x,m)dm\nonumber \\
u(t=0,x,m) & = & \phi(x,m)\nonumber \\
u(t,x,m=1) & = & 0\nonumber \\
w(t,x=-\infty)=1 &  & w(t,x=+\infty)=0.\nonumber 
\end{eqnarray}
In Section \ref{sec:Numerical-simulations}, we will consider the
full model \eqref{eq:structured_fisher_eqn} when the rates of cellular
diffusion and proliferation are piece-wise constant functions of $m$
and numerically investigate how different functions for $f(t)$ lead
to increased and decreased levels of population migration.

\section{Background Material from Size-Structured Population Modeling\label{sec:Characteristic-equations-from}}

Before investigating the existence of traveling-wave solutions to
\eqref{eq:structured_fisher_eqn_norm}, it is useful to review some
key topics used to solve size-structured population models, as discussed
in \citep{webb_population_2008}. These topics will be useful in analyzing
\eqref{eq:structured_fisher_eqn} in later sections. A reader familiar
with using the method of characteristics to solve size-structured
population models may briefly skim over this section to pick up on
the notation used throughout our study.

As an example, we consider the size-structured model given by
\begin{eqnarray}
u_{t}+(g(y)u)_{y} & = & Au\label{eq:size_example}\\
u(t=0,y) & = & \phi(y)\nonumber 
\end{eqnarray}
where $u=u(t,y):[0,\infty)\times[y_{0},y_{1}]\rightarrow\mathbb{R}$
denotes the size distribution over $y$ of a population at time $t$,
$y_{0}$ and $y_{1}$ denote the minimum and maximum population sizes
respectively, and $g(y)\in C^{1}(([y,y_{1}])$ denotes the physical
growth rate\footnote{Note that in this section, $g(y)$ denotes a growth rate with respect
to size, $y$, whereas throughout the rest of our study, $g(m)$ denotes
an activation rate with respect to biochemical activity, $m$.} of individuals of size $y$. In this section, we will work in the
Banach space $\mathbb{X}=L^{1}((y_{0},y_{1})\rightarrow\mathbb{R}),$
and assume $A\in\mathbb{\mathcal{B}}(\mathbb{X})$, the space of bounded,
linear operators on $\mathbb{X}$. The method of characteristics will
facilitate solving \eqref{eq:size_example}.

For a fixed size $\underline{y}\in[y_{0},y_{1}],$ the function
\begin{equation}
\sigma(y;\underline{y}):=\int_{\underline{y}}^{y}\frac{1}{g(y')}dy'\label{eq:sigma}
\end{equation}
provides \emph{the time it takes for an individual to grow from the
fixed size $\underline{y}$ to arbitrary size $y$.} If $g(y)$ is
positive and uniformly continuous on $[y_{0},y_{1}]$, then $\sigma(y;\underline{y})$
is invertible. We denote the inverse function, $\sigma^{-1}(t;\underline{y}),$
as the \emph{growth curve}, and it computes \emph{the size of an individual
over time that starts at size $\underline{y}$ at time $t=0$}. For
instance, if an individual has size $\underline{y}$ at $t=0,$ then
that individual will have size $\sigma^{-1}(t_{1};\underline{y})$
at time $t=t_{1}$. Some helpful properties of the growth curve are
that $\sigma^{-1}(0;\underline{y})=\underline{y}$ and 
\begin{equation}
\frac{d}{dt}\sigma^{-1}(t;\underline{y})=g(\sigma^{-1}(t;\underline{y})).\label{eq:dsigma/dt}
\end{equation}
See Section \ref{sec:Properties-of} in the appendix for the derivation
of \eqref{eq:dsigma/dt}. 

In order to solve \eqref{eq:size_example} with the method of characteristics,
we set $y=\sigma^{-1}(t;\underline{y})$ to define the variable $v(t;\underline{y})$:
\begin{equation}
v(t;\underline{y}):=u(t,y=\sigma^{-1}(t;\underline{y})).\label{eq:size_characteristic}
\end{equation}
As shown in Section \ref{sec:Derivation-of} of the appendix, substitution
of \eqref{eq:size_characteristic} into \eqref{eq:size_example} yields
the characteristic equation
\begin{equation}
v_{t}=-g'(\sigma^{-1}(t;\underline{y}))v+Av,\label{eq:v_t_diffeq}
\end{equation}
where primes denote differentiation with respect to $y$. This characteristic
equation has size $\underline{y}$ at time $t=0$ and can be solved
explicitly as\footnote{To derive this, use separation of variables and with the help of \eqref{eq:dsigma/dt}
note that $\int_{0}^{t}g'(\sigma^{-1}(\tau;\underline{y}))d\tau=\ln[g(\sigma^{-1}(t;\underline{y}))/g(\underline{y})]$.} 
\begin{equation}
v(t;\underline{y})=\frac{g(\underline{y})}{g(\sigma^{-1}(t;\underline{y}))}e^{At}\phi(\underline{y}).\label{eq:soln_size_characteristic}
\end{equation}
As \eqref{eq:soln_size_characteristic} provides the solution to \eqref{eq:size_example}
along the arbitrary characteristic curve with initial size $\underline{y},$
we use it to solve the whole equation with the substitution $y=\sigma^{-1}(t,\underline{y})$,
in which we find
\begin{equation}
u(t,y)=\left\{ \begin{array}{cc}
\frac{g(\sigma^{-1}(-t,y))}{g(y)}e^{At}\phi(\sigma^{-1}(-t,y)) & \sigma^{-1}(t;y_{0})\le y\le y_{1}\\
0 & y_{0}\le y<\sigma^{-1}(t;y_{0}).
\end{array}\right.\label{eq:soln_to_size_example}
\end{equation}
If $\phi(y)\notin C^{1}(y_{0},y_{1}),$ then \eqref{eq:soln_to_size_example}
is viewed as a generalized solution. Note that a piecewise form is
needed for \eqref{eq:soln_to_size_example} because we do not have
any individuals below the minimum size, $y_{0},$ and thus the minimum
possible size at time $t$ is given by $\sigma^{-1}(t;y_{0}).$ If
the population is assumed to give birth to individuals of size $y_{0}$
over time, then the appropriate renewal equation representing population
birth would replace the zero term in the piecewise function (see \citep[$\S$ 9.5]{banks_mathematical_2009}
for an example in size-structured populations and \citep{gourley_vector_2007}
for an example in age-structured populations).

\section{Existence of Traveling Wave Solutions to the Structured Fisher's
Equation\label{sec:Traveling-waves-analysis}}

\subsection{Existence of traveling wave solutions to \eqref{eq:structured_fisher_eqn_norm}}

We now incorporate topics from the previous section to show the existence
of traveling wave solutions to \eqref{eq:structured_fisher_eqn_norm}.
 After taking the time derivative of $w(t,x)$, which was defined
in \eqref{eq:w_def}, we can rewrite \eqref{eq:structured_fisher_eqn_norm}
as a system of two coupled PDEs\footnote{Note that either $g(m_{0})=0$ or $u(t,m=m_{0},x)=0$ for $t>0$,
so that the activation term drops out when integrating over $m$ for
$w$.}:
\begin{eqnarray}
u_{t}+(g(m)u)_{m} & = & u_{xx}+u(1-w)\nonumber \\
w_{t} & = & w_{xx}+w(1-w).\label{eq:coupled_fisher_uw}
\end{eqnarray}
Note that in this section, $g(m)$ is a function of biochemical activity
level and $\sigma^{-1}(t;\underline{m})$ computes \emph{the activity
level of an individual over time that starts at level $\underline{m}$
at time $t=0$}. We will thus now refer to $\sigma^{-1}(t;\underline{m})$
as the \emph{activation curve}. We next set up the characteristic
equation for $u$ by setting $m=\sigma^{-1}(t;\underline{m})$ for
a fixed value of $\underline{m}$:
\begin{equation}
v(t,x;\underline{m}):=u(t,x,m=\sigma^{-1}(t;\underline{m})).\label{eq:v_characteristic}
\end{equation}
Substituting \eqref{eq:v_characteristic} into \eqref{eq:coupled_fisher_uw}
simplifies to our characteristic equation
\begin{eqnarray}
v_{t} & = & v_{xx}+v[1-w-g'(\sigma^{-1}(t;\underline{m}))]\nonumber \\
w_{t} & = & w_{xx}+w(1-w),\label{eq:coupled_fisher_vw}
\end{eqnarray}
a nonautonomous system of two coupled PDEs in time and space. Note
that the bottom equation for \eqref{eq:coupled_fisher_vw} is Fisher's
Equation, which has positive monotonic traveling wave solutions for
any speed $c\ge2$ (see \citep[$\S$ 11.2]{murray_mathematical_2002}). 

We next aim to derive traveling wave solutions to \eqref{eq:coupled_fisher_vw},
however, we are not aware of any traveling wave solutions to nonautonomous
systems such as this one. From our knowledge of size-structured population
models from Section \ref{sec:Characteristic-equations-from}, we instead
intuit the ansatz of a self-similar traveling wave solution, which
we write as 
\begin{eqnarray}
v(t,x;\underline{m}) & = & \frac{g(\underline{m})}{g(\sigma^{-1}(t;\underline{m}))}V(z),\ z=x-ct\label{eq:self_similar_TW_ansatz}\\
w(t,x) & = & W(z).\nonumber 
\end{eqnarray}
In this ansatz, $V(z)$ will define a traveling wave profile for $v$
and $\frac{g(\underline{m})}{g(\sigma^{-1}(t;\underline{m}))}$ will
provide the height of the function over time. With the aid of the
chain rule, we observe that:
\begin{eqnarray*}
v_{t}(t,x;\underline{m}) & = & \frac{g'(\sigma^{-1}(t;\underline{m}))g(\sigma^{-1}(t;\underline{m}))}{g(\underline{m})}V-c\frac{g(\sigma^{-1}(t;\underline{m}))}{g(\underline{m})}V_{z}\\
v_{xx}(t,x;\underline{m}) & = & \frac{g(\sigma^{-1}(t;\underline{m}))}{g(\underline{m})}V_{zz},
\end{eqnarray*}
where subscripts denote differentiation with respect to t, $x$, or
$z$ and primes denote differentiation with respect to $m$. Substituting
\eqref{eq:self_similar_TW_ansatz} into \eqref{eq:coupled_fisher_vw}
reduces to the autonomous system 
\begin{eqnarray}
-cV_{z} & = & V_{zz}+V(1-W)\nonumber \\
-cW_{z} & = & W_{zz}+W(1-W).\label{eq:coupled_fisher_VW}
\end{eqnarray}

It is now useful to rewrite \eqref{eq:coupled_fisher_VW} as the first
order system

\begin{equation}
\frac{d}{dz}\boldsymbol{\mathcal{\mathcal{V}}}=\left(\begin{array}{c}
V_{z}\\
-cV_{z}-V(1-W)\\
W_{z}\\
-cW_{z}-W(1-W)
\end{array}\right)\label{eq:first_order_VW}
\end{equation}
for $\boldsymbol{\mathcal{V}}(z)=[V(z),V_{z}(z),W(z),W_{z}(z)]^{T}$.
Recall that profiles to traveling wave solutions can be constructed
with heteroclinic orbits between equilibria for a given dynamical
system (or homoclinic orbits for a traveling pulse) \citep[$\S$ 6.2]{keener_mathematical_2009}.
We observe two types of equilibria for \eqref{eq:coupled_fisher_VW},
given by $\boldsymbol{\mathcal{\mathcal{V}}}_{1}^{*}=(1,0,V,0)^{T}$
and $\boldsymbol{\mathcal{\mathcal{V}}}_{2}^{*}=\vec{0}$, where the
former represents a confluent cell density and the latter represents
an empty wound space. We accordingly search for heteroclinic orbits
from $\boldsymbol{\mathcal{\mathcal{V}}}_{1}^{*}$ to $\boldsymbol{\mathcal{\mathcal{V}}}_{2}^{*}$
for some $c>0.$ We choose to focus on the characteristic equations
$v(t,x,m=\sigma^{-1}(t;\underline{m}))$ for values of $\underline{m}$
in which $\phi(\underline{m},x=-\infty)>0$ to represent the population
of cells migrating into the empty wound space. We thus denote $\boldsymbol{\mathcal{\mathcal{V}}}_{1}^{*}=(1,0,\nu,0)^{T}$
for $\nu>0$. 

Note that $\boldsymbol{\mathcal{\mathcal{V}}}_{1}^{*}$ is an equilibrium
for any value of $V,$ as $W=1$ will guarantee the existence of an
equilibrium. Such a ``continuum'' of equilibria was also observed
in \citep{perumpanani_traveling_2000}. This structure of $\boldsymbol{\mathcal{\mathcal{V}}}_{1}^{*}$
yields a zero eigenvalue after linearizing \eqref{eq:first_order_VW}
about $\boldsymbol{\mathcal{\mathcal{V}}}_{1}^{*}$, so we cannot
use linear theory to study the local behavior of \eqref{eq:first_order_VW}
near $\boldsymbol{\mathcal{\mathcal{V}}}_{1}^{*}$. While we could
construct the unstable manifold of \eqref{eq:first_order_VW} using
a power series representation to study its local behavior around $\boldsymbol{\mathcal{\mathcal{V}}}_{1}^{*}$
(see \citep[Section 5.6]{meiss_differential_2007}), we find it more
insightful to define a trapping region in the $(V,V_{z})$-plane as
has been done in previous traveling wave studies \citep{ai_reaction_2015,kuang_data-motivated_2015}.
We will then use asymptotically autonomous phase-plane theory to describe
the $\omega$-limit set of our flow, which will show the existence
of a heteroclinic orbit from $\boldsymbol{\mathcal{\mathcal{V}}}_{1}^{*}$
to $\boldsymbol{\mathcal{\mathcal{V}}}_{2}^{*}$. \emph{Trapping regions}
are positively invariant regions with respect to the flow of a dynamical
system, and the \emph{$\omega$-limit set }of a flow is the collection
all limit point of that flow \citep[$\S$ 4.9-10]{meiss_differential_2007}.

We study the trajectory of $\boldsymbol{\mathcal{\mathcal{V}}}$ in
the $(V,V_{z})$-plane by defining the triangular region bound by
the lines $\{V=\nu,V_{z}=0,V_{z}=-\frac{c}{2}V\}$ and denoting this
region as $\Delta.$ The following lemma will demonstrate that $\Delta$
is a trapping region for the flow of \eqref{eq:first_order_VW} in
the $(V,V_{z})$-plane. 

\textbf{Lemma: }Let $\nu>0$ and $c\ge2$. Then the region $\Delta$
is positively invariant with respect to \eqref{eq:first_order_VW}
so long as $0<W(z)<1$ for all $z\in\mathbb{R}.$

\emph{Proof:}

We prove this lemma by investigating the vector field along each of
the lines specifying the boundary of our region and showing that they
point into the interior of the space.

i.) Along the line $V_{z}=0,$ $\frac{d}{dz}V_{z}=-V(1-W)$, which
is nonpositive because $W(z)<1$ for all $z\in\mathbb{R}$ and our
region is defined for $V(z)\ge0.$ If $V=0,$then we are at the equilibrium
point $(V,V_{z})=(0,0)$.

ii.) Along $V=\nu,$ $\frac{d}{dz}V=V_{z},$ which is negative in
our defined region. The only point to worry about here is at $(V,V_{z})=(\nu,0)$,
as then $\frac{d}{dz}V=0$. However, we see from part i.) that $\frac{d}{dz}V_{z}<0$
here, so that a flow starting at $(\nu,0)$ will initially move perpendicular
to the $V$-axis in the negative $V_{z}$ direction, and then $\frac{d}{dz}V<0$,
so the flow enters $\Delta.$

iii.) Note that the inner normal vector to the line $V_{z}=-\frac{c}{2}V$
is $\hat{n}=\left(\frac{c}{2},1\right).$ Then 
\begin{eqnarray*}
\hat{n}\cdot\frac{d}{dz}(V,V_{z}) & = & \left(\frac{c}{2},1\right)\cdot(V_{z},-cV_{z}-V+VW)\\
 & = & \left(\frac{c}{2},1\right)\cdot(-\frac{c}{2}V,\frac{c^{2}}{2}V-V+VW)\\
 & = & -\frac{c^{2}}{4}V+\frac{c^{2}}{2}V-V+VW\\
 & = & V\left(\frac{c^{2}}{4}-1\right)+VW,
\end{eqnarray*}
which is positive, as $c\ge2$.$\square$

\begin{figure}
\includegraphics[width=0.45\textwidth]{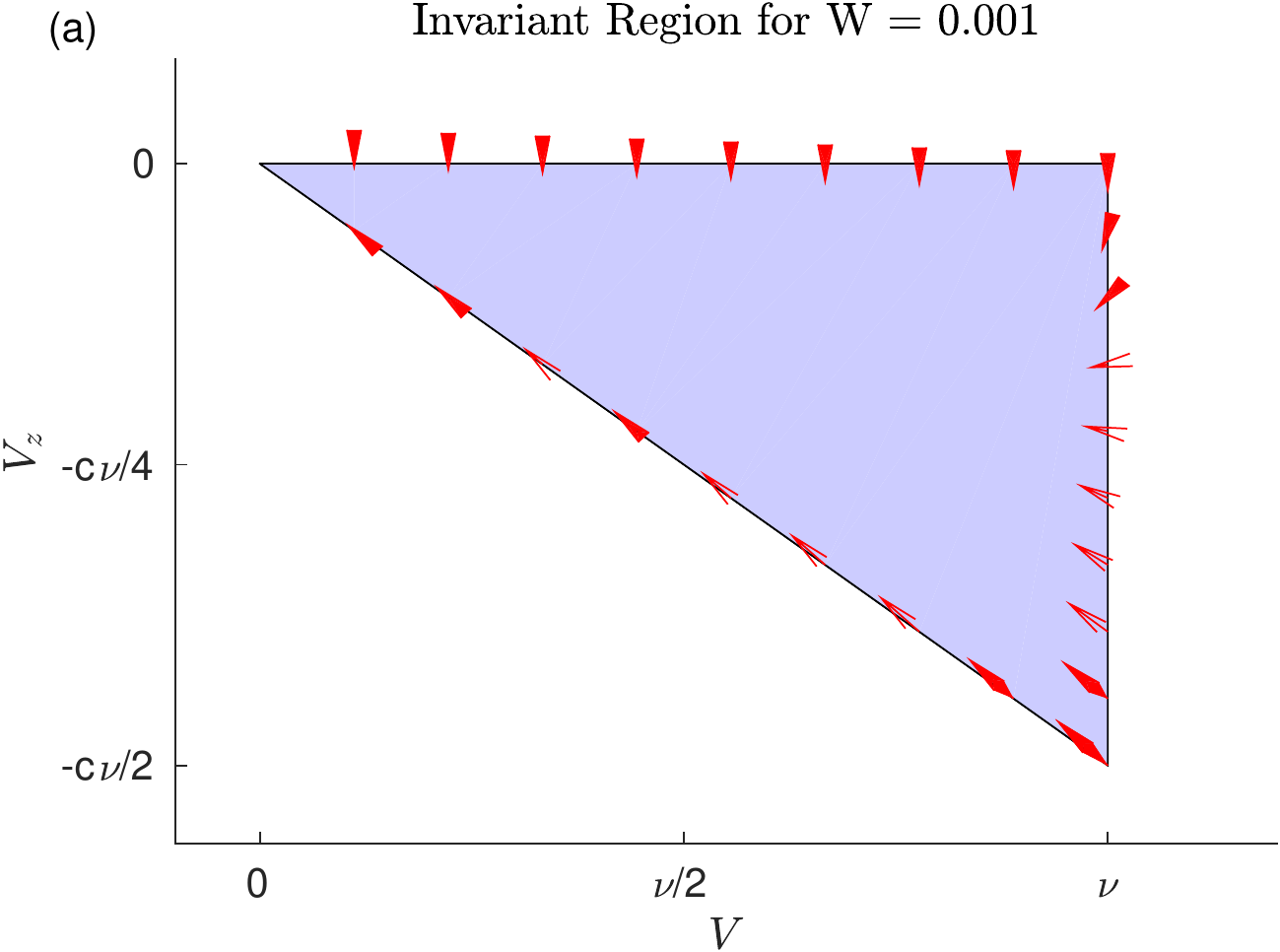}\hfill{}\includegraphics[width=0.45\textwidth]{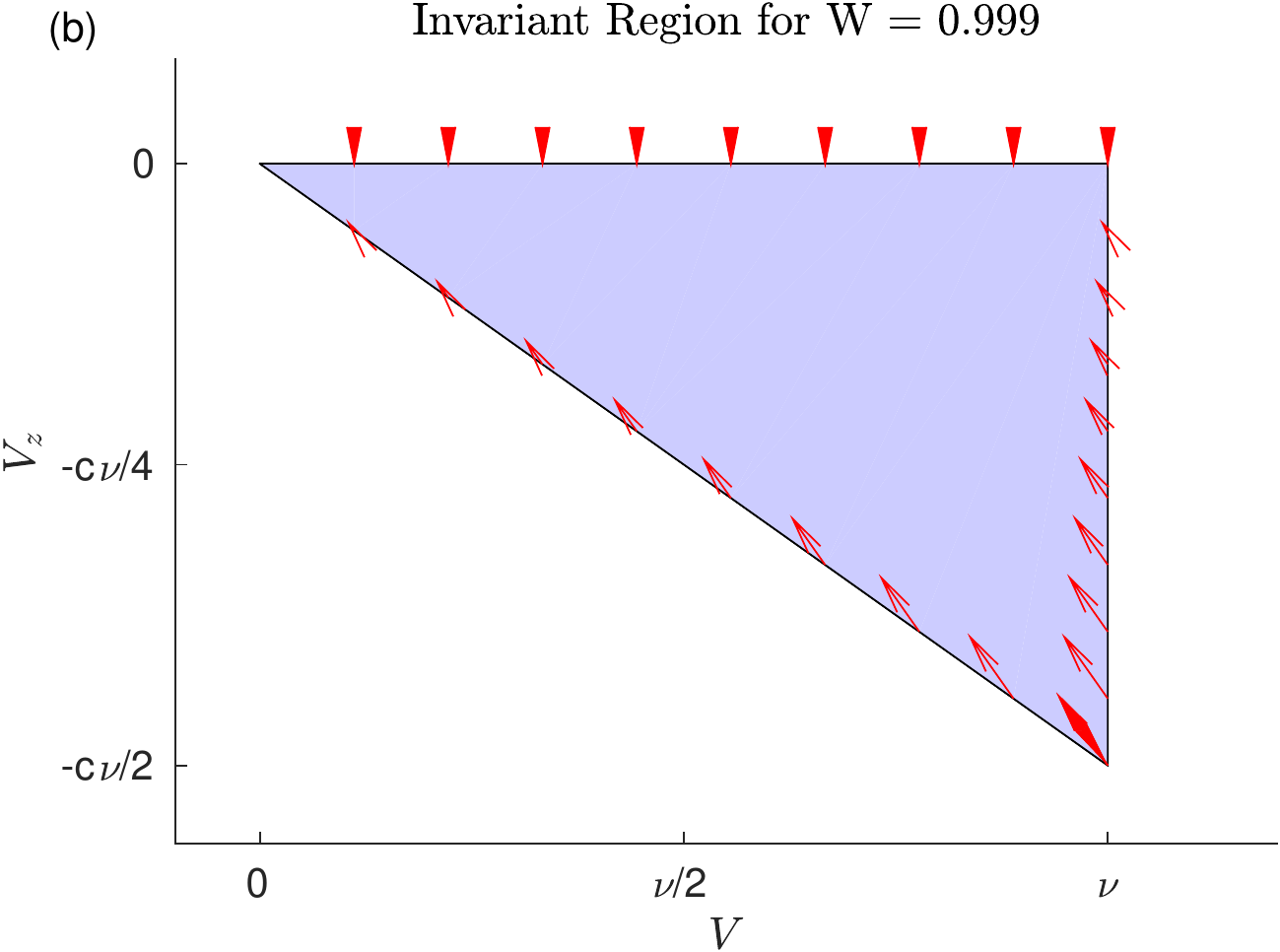}\vspace{.1cm}

\includegraphics[width=0.45\textwidth]{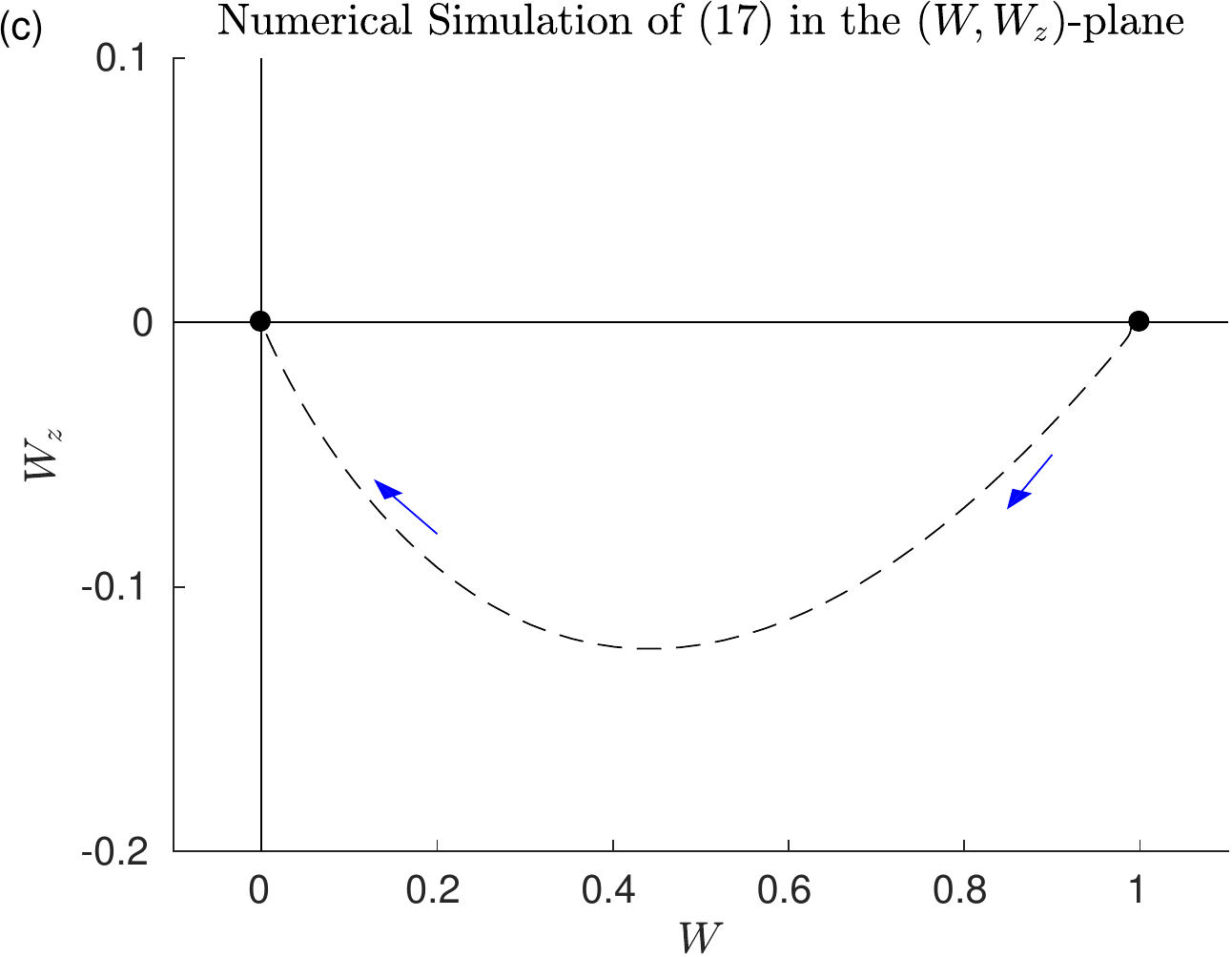}\hfill{}\includegraphics[width=0.45\textwidth]{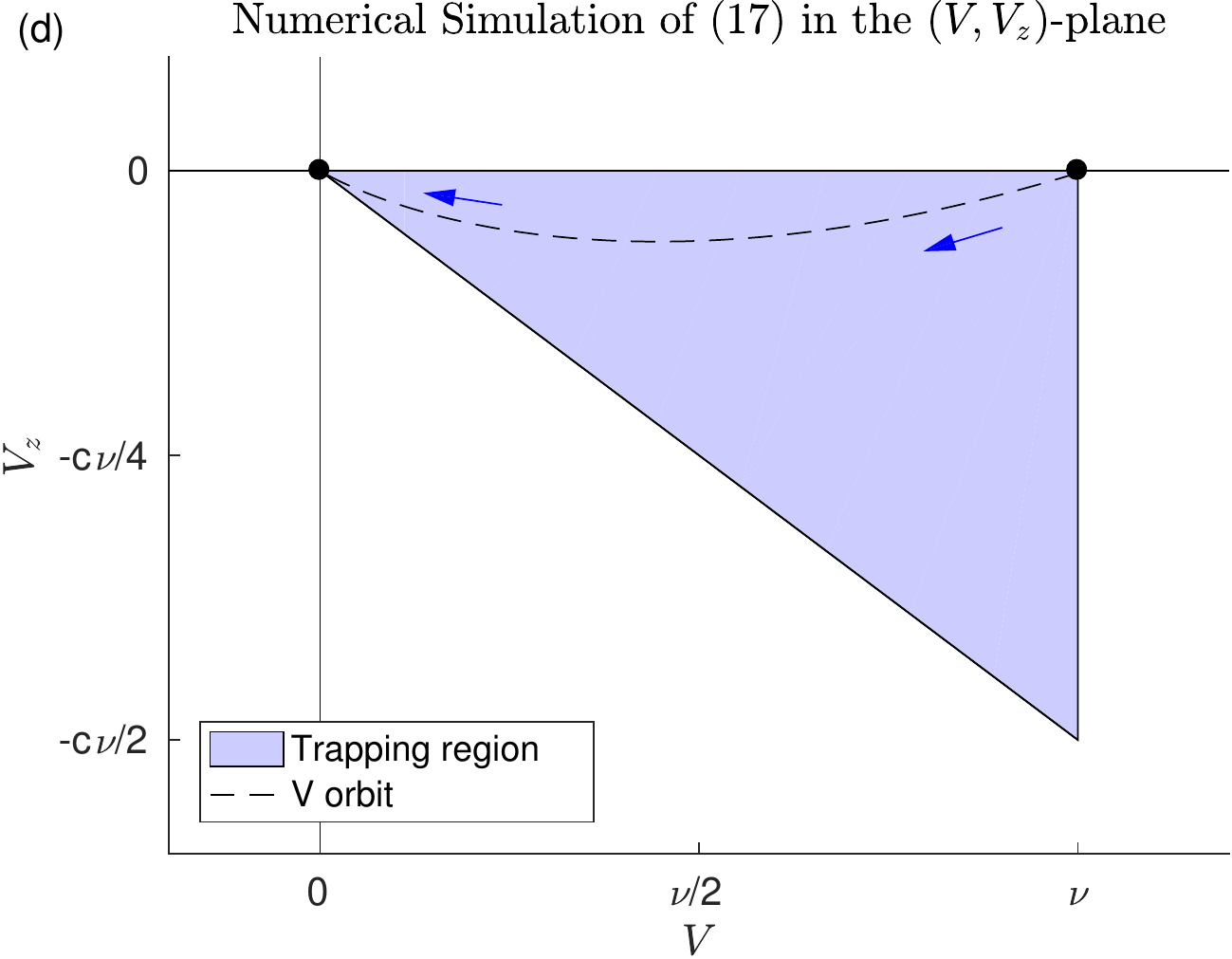}

\protect\caption{Construction of the heteroclinic orbit between $\boldsymbol{\mathbb{\mathcal{V}}}_{1}^{*}$
and $\boldsymbol{\mathbb{\mathcal{V}}}_{2}^{*}$ for \eqref{eq:first_order_VW}.
In (a) and (b), we depict the trapping region in the $(V,V_{z})$-plane,
$\Delta$, and the vector field along its boundary for $\alpha=0.5,c=2.3,W_{z}=0$
and $W$ near 0 and 1, respectively. In (c) and (d), we depict numerical
simulations of \eqref{eq:first_order_VW} in the $(W,W_{z})$-plane
and in the $(V,V_{z})$-plane, respectively. Arrows denote the direction
of the flow and the black dots mark equilibria of \eqref{eq:first_order_VW}.
\label{fig:invariant_region}}
\end{figure}

This proof is visually demonstrated in the top row of Figure \ref{fig:invariant_region}.
As $W(z)$ has a heteroclinic orbit with $W(-\infty)=1$ and $W(\infty)=0$
for any $c\ge2$ \citep[$\S$ 11.2]{murray_mathematical_2002}, we
conclude that $\Delta$ is a positively invariant set for the flow
of \eqref{eq:first_order_VW} in the $(V,V_{z})$. The following corollary
describes the $\omega$-limit set of \eqref{eq:first_order_VW}.

\textbf{Corollary: }The $\omega$-limit set of \eqref{eq:first_order_VW}
starting at $\boldsymbol{\mathbb{\mathcal{V}}}_{2}^{*}$ , $\omega(\boldsymbol{\mathbb{\mathcal{V}}}_{2}^{*})$,
is $\boldsymbol{\mathbb{\mathcal{V}}}_{1}^{*}$$.$

\emph{Proof:}

As $W(z)\rightarrow0$ as $z\rightarrow+\infty,$ then the vector
field for \eqref{eq:first_order_VW} in the $(V,V_{z})$-plane is
asymptotically autonomous to the vector field
\begin{equation}
\frac{d}{dz}\left(\begin{array}{c}
V\\
V_{z}
\end{array}\right)=\left(\begin{array}{c}
V_{z}\\
-cV_{z}-V
\end{array}\right),\label{eq:asymptotical_flow}
\end{equation}
a linear system whose only equilibrium is the origin. As $c\ge2,$
the origin is a stable equilibrium and the flow of the limiting system
remains in $\Delta$, and hence the fourth quadrant, for all time.\\
As $\nicefrac{d}{dz}V=V_{z}<0$ in $\Delta$, no periodic or homoclinic
orbits can exist for the limiting system. We thus conclude from the
asymptotically autonomous Poincare-Bendixson Theorem presented in
\citep{markus_asymptotically_1956}\footnote{The relevant theorem statement is given in Appendix \ref{sec:Statement-of-AA_markus}.
We note that while the results of \citep{markus_asymptotically_1956}
are sufficient for our study, asymptotically autonomous systems have
been more extensively studied in \citep{blythe_autonomous_1991,castillo-chavez_asymptotically_1994,thieme_convergence_1992,thieme_asymptotically_1993}
and a more comprehensive result in describing the $\omega$-limit
set of the asymptotically autonomous flow is given in \citep{thieme_asymptotically_1993}.} that our flow in the $(V,V_{z})$ plane starting at $(\nu,0)$ will
limit to the origin. We conclude that $\omega(\boldsymbol{\mathbb{\mathcal{V}}}_{2}^{*})=\boldsymbol{\mathbb{\mathcal{V}}}_{1}^{*}.$$\square$

\subsection{Summary of results}

We have demonstrated the existence of self-similar traveling wave
solutions to \eqref{eq:structured_fisher_eqn_norm} in this section
of the form
\begin{eqnarray*}
u(t,x,m) & = & \frac{g(\underline{m})}{g(\sigma^{-1}(t;\underline{m}))}V(z;\underline{m}),\ z=x-ct
\end{eqnarray*}
for $c\ge2$ and values of $\underline{m}$ where the initial condition
$\phi(x=-\infty,\underline{m})>0$ and $V(z;\underline{m})=u(x-ct,m=\sigma^{-1}(t;\underline{m})).$
Setting $m=\sigma^{-1}(t,\underline{m}),$ then this can be written
more explicitly as
\[
u(t,x,m)=\left\{ \begin{array}{cc}
\frac{g(\sigma^{-1}(-t,m))}{g(m)}V(z;\sigma^{-1}(-t;m)), & \sigma^{-1}(t;0)\le m<1\\
0 & \mbox{otherwise}
\end{array}\right.
\]
and 
\[
\int_{0}^{1}u(t,x,m)dm=w(t,x)=W(z)
\]
where $[V(z;\underline{m}),W(z)]^{T}$ satisfies \eqref{eq:coupled_fisher_VW}.
An example height function $\frac{g(\sigma^{-1}(-t,m))}{g(m)}$ for
trajectories along the activation curves $m=\sigma^{-1}(t;\underline{m})$
will be demonstrated later in Figure \ref{fig:p(t,m)_ex1}.

\section{Structured Fisher's Equation with MAPK-dependent Phenotype\label{sec:Numerical-simulations}}

We now study a version of Fisher's Equation where cellular migration
and proliferation depend on biochemical activity, $m$. Various cell
lines have reduced rates of proliferation and increased migration
in response to MAPK activation \citep{chapnick_leader_2014,clark_molecular_1995,matsubayashi_erk_2004},
so we let $m$ denote activity along the MAPK signaling cascade in
this section. We consider a model with two subpopulations: one with
a high rate of diffusion in response to MAPK activation and the other
with a high rate of proliferation when MAPK levels are low. MAPK activation
will depend on an external forcing factor to represent the presence
of an extracellular signaling chemical, such as ROS, TGF-$\beta$,
or EGF. While the method of characteristics is not applicable to spatial
activation patterning here due to the parabolic nature of \eqref{eq:structured_fisher_eqn}
in space, we can investigate temporal patterns of activation and deactivation.
\emph{We will exhibit simple scenarios that give rise to three ubiquitous
patterns of biochemical activity: 1.) a sustained wave of activation,
2.) a single pulse of activitation, and 3.) periodic pulses of activation.}

Before describing these examples, we first introduce some tools to
facilitate our study of \eqref{eq:structured_fisher_eqn}. We will
detail some assumptions that simplify our analysis in Section \ref{sub:Model-Description},
solve and compute the population activation profile over time and
use it to define some activation criteria in Section \ref{sub:Activation-level-in},
and discuss numerical issues and the derivation of a nonautonomous
averaged Fisher's Equation in Section \ref{sub:Derivation-of-the}
before illustrating the different activation patterns and their effects
on migration in Section \ref{sub:Three-biologically-motivated-exa}.

\subsection{Model Description\label{sub:Model-Description}}

Recall that the full structured Fisher's Equation is given by

\begin{eqnarray}
u_{t}+(f(t)g(m)u)_{m} & = & D(m)u_{xx}+\lambda(m)u\left(1-w\right)\label{eq:nonautonomous_fisher_eqn}\\
w & = & \int_{m_{0}}^{m_{1}}u(t,x,m)dm\nonumber \\
u(t=0,x,m) & = & \phi_{1}(m)\phi_{2}(x)\nonumber \\
u(t,x,m=m_{1}) & = & 0\nonumber \\
w(t,x=+\infty)=0 &  & w(t,x=-\infty)=1.\nonumber 
\end{eqnarray}
We have chosen the separable initial condition $u(t=0,x,m)=\phi_{1}(m)\phi_{2}(x)$
for simplicity. Given some $m_{crit}\in(m_{0},m_{1})$, we define
two subsets of $[m_{0},m_{1}]$ as $M_{inact}:=[m_{0},m_{crit}],M_{act}:=(m_{crit},m_{1}]$,
and the rates of diffusion and proliferation by

\begin{equation}
D(m):=\left\{ \begin{array}{cc}
D_{1} & m\in M_{inact}\\
D_{2} & m\in M_{act}
\end{array}\right.,\ \lambda(m):=\left\{ \begin{array}{cc}
\lambda_{1} & m\in M_{inact}\\
\lambda_{2} & m\in M_{act}
\end{array}\right.\label{eq:nonatuonomous_D_lambda}
\end{equation}
 for $D_{1}<D_{2}$ and $\lambda_{2}<\lambda_{1}.$ Hence for $m\in M_{inact},$
the population is termed as \emph{inactive }and primarily proliferates
whereas for $m\in M_{act},$ the population is termed as \emph{active}
and primarily diffuses.

We let $supp(\phi_{1}(m))=[\underline{m}_{\min},\mbox{\ensuremath{\underline{m}}}_{\max}]$
for $\underline{m}_{\max}<m_{crit}$ and assume that $\int_{m_{0}}^{m_{1}}\phi_{1}(m)dm=1$
so $\phi_{1}(m)$ represents a probability density function for the
initial distribution of cells in $m$. We accordingly denote 
\[
\Phi_{1}(m):=\left\{ \begin{array}{cc}
0 & m\le m_{0}\\
\int_{m_{0}}^{m}\phi_{1}(m')dm' & m_{0}<m\le m_{1}\\
1 & m_{1}<m
\end{array}\right.
\]
as the cumulative distribution function for $\phi_{1}(m).$

\subsection{Activation Profile and Activation Criteria\label{sub:Activation-level-in}}

An interesting question is how the distribution of \eqref{eq:nonautonomous_fisher_eqn}
along $m$ changes over time. To answer this question, we consider
\eqref{eq:nonautonomous_fisher_eqn} in terms of $t$ and $m$, which
we will write as $p(t,m)$ and call the \emph{activation profile}:
\begin{eqnarray}
p_{t}+(f(t)g(m)p)_{m} & = & 0\label{eq:dist_p_tm}\\
p(0,m) & = & \phi_{1}(m).\nonumber 
\end{eqnarray}
Following the analysis from Section \ref{sec:Characteristic-equations-from},
we can solve \eqref{eq:dist_p_tm} analytically. We integrate \eqref{eq:dist_p_tm}
along the \emph{activation curve}s, $h(t;\underline{m}),$ which now
are given by 
\begin{equation}
m=h(t;\underline{m}):=\sigma^{-1}\left(F(t);\underline{m}\right),\label{eq:h_characteristic}
\end{equation}
where $F(t):=\int_{0}^{t}f(\tau)d\tau$ denotes a cumulative activation
function. We find the activation profile to be:
\begin{equation}
p(t,m)=\left\{ \begin{array}{cc}
\frac{g(\sigma^{-1}(-F(t),m))}{g(m)}\phi_{1}(\sigma^{-1}(-F(t),m)) & h(t;m_{0})\le m\le m_{1}\\
0 & m_{0}\le m<h(t;m_{0}).
\end{array}\right.\label{eq:u(tm)_soln_exact_fs}
\end{equation}

Now we can derive a condition for a cell population starting in the
inactive population to enter the active population. We see from \eqref{eq:u(tm)_soln_exact_fs}
that the population will enter the active population if 
\begin{equation}
h(t;\underline{m}_{\max})>m_{crit}\iff F(t)>\sigma(m_{crit};\underline{m}_{\max})\label{eq:transition_condition}
\end{equation}
 for some values of $t$. By standard calculus arguments, \eqref{eq:transition_condition}
will occur if
\begin{equation}
F(t_{\max})>\sigma(m_{crit};\underline{m}_{\max})\label{eq:transition_criterion}
\end{equation}
where a local maximum for $F(t)$ occurs at $t=t_{\max}$. Hence,
$f(t_{\max})=0$, $f(t_{\max}^{-})>0$, and $f(t_{\max}^{+})<0$.
We denote \eqref{eq:transition_criterion} as the \emph{activation
criterion} for \eqref{eq:nonautonomous_fisher_eqn}. By the same argument,
for the entire population to activate at some point, then we can derive
the \emph{entire activation criterion} as
\begin{equation}
F(t_{\max})>\sigma(m_{crit};\underline{m}_{\min}).\label{eq:full_transition_criterion}
\end{equation}

\subsection{Numerical Simulation Issues and Derivation of an Averaged Nonautonomous
Fisher's Equation\label{sub:Derivation-of-the}}

\begin{figure}
\subfloat[]{\protect\includegraphics[width=0.45\textwidth]{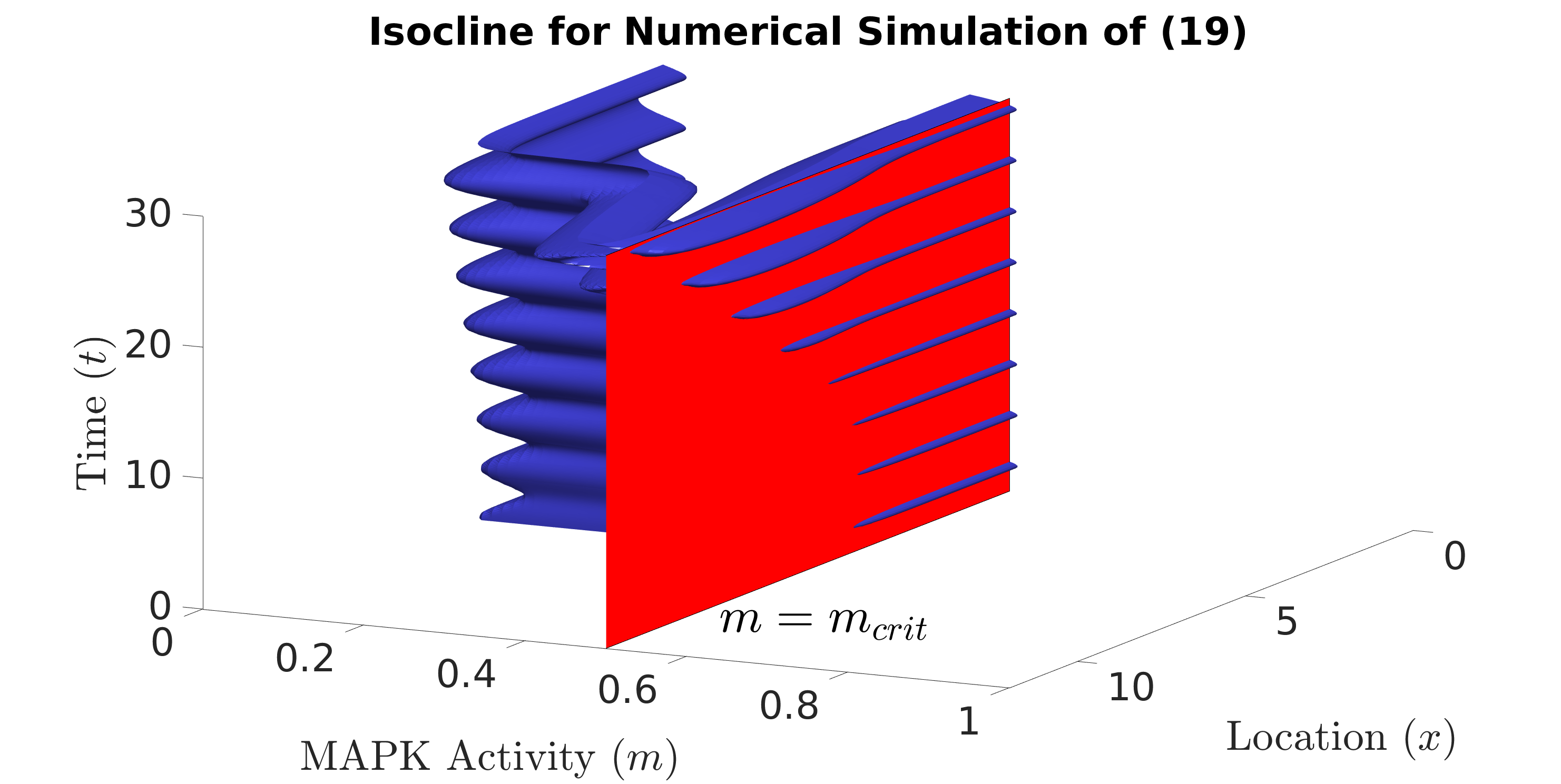}}\subfloat[]{\protect\includegraphics[width=0.45\textwidth]{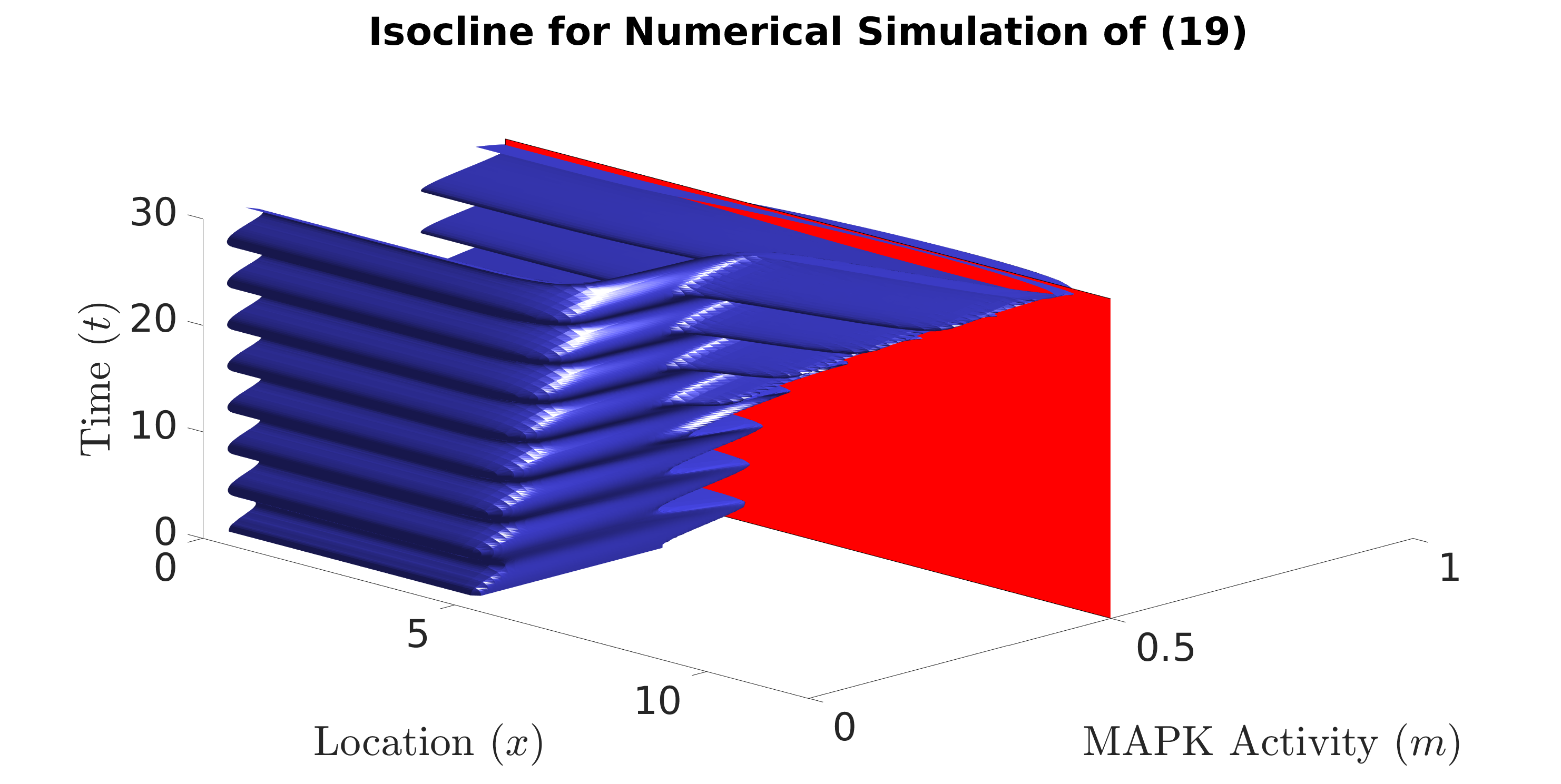}}

\centering{}\protect\caption{Two views of the isocline for $u=1$ from a numerical simulation of
\eqref{eq:nonautonomous_fisher_eqn} with $g(m)=\alpha m(1-m)$ and
$f(t)=\beta\sin(\gamma t)$ for $\alpha=0.5,\beta=1,$ $\gamma=1.615,D_{1}=0.01,D_{2}=1,\lambda_{1}=0.25,$
and $\lambda_{2}=0.0025$ and an initial condition of $\phi_{1}(m)=\nicefrac{10}{3}I_{[.05,0.35]}(m)$
and $\phi_{2}(x)=I_{[x\le5]}(x)$. The numerical scheme is discussed
in Section \ref{sub:Derivation-of-the} and the step sizes used are
$\Delta m=\nicefrac{1}{80},\Delta x=\nicefrac{1}{5},\Delta t=10^{-3}.$
From \eqref{eq:transition_criterion}, the simulation should not cross
the $m=m_{crit}$ plane, which is given by the red plane. We see in
frame (a) that the simulation does cross the $m=m_{crit}$ plane due
to numerical diffusion, which causes the high rate of diffusion along
$x$ seen in frame (b). \label{fig:Numerical-simulations-foru(tmx)}}
\end{figure}

We depict the $u=1$ isocline for a numerical simulation of \eqref{eq:nonautonomous_fisher_eqn}
in Figure \ref{fig:Numerical-simulations-foru(tmx)} with $g(m)=\alpha m(1-m)$,
$f(t)=\beta\sin(\gamma t)$, $\alpha=1/2,\beta=1,$ and $\gamma=1.615$.
These terms will be detailed more in Example 3 below. For numerical
implementation, we use a standard central difference scheme for numerical
integration along the $x$-dimension, an upwind scheme with flux limiters
(similar to those described in \citep{thackham_computational_2008})
to integrate along the $m$ dimension, and a Crank-Nicholson scheme
to integrate along time. From \eqref{eq:transition_criterion}, we
see that this simulation should not enter the active population with
an initial condition of $\phi_{1}(m)=\nicefrac{10}{3}I_{[.05,0.35]}(m),$
where $I_{M}(m)$ denotes an indicator function with support for $m\in M$.
In Figure \ref{fig:Numerical-simulations-foru(tmx)}, however, we
observe that the numerical simulation does enter the active population,
which causes a significant portion of the population to incorrectly
diffuse into the wound at a high rate.

Numerical simulations of advection-driven processes have been described
as an ``embarrassingly difficult'' task, and one such problem is
the presence of numerical diffusion \citep{leonard_ultimate_1991,thackham_computational_2008}.
Numerical diffusion along the $m$-dimension is hard to avoid and
here causes a portion of the cell population to enter the active population
in situations where the it should approach the $m=m_{crit}$ plane
but not pass it. Numerical diffusion can be reduced with a finer grid,
but this can lead to excessively long computation times. With the
aid of the activation curves given by \eqref{eq:h_characteristic},
however, we can track progression of cells in the $m$-dimension analytically
and avoid the problems caused by numerical diffusion completely. 

To avoid the problems caused by numerical diffusion, we derive a nonautonomous
Fisher's Equation for $w(t,x)$ that represents the average behavior
along $m$ with time-dependent diffusion and proliferation terms.
To investigate the averaged cell population behavior along $m$ over
time, we integrate \eqref{eq:nonautonomous_fisher_eqn} over $m$
to find
\begin{eqnarray}
w_{t}(t,x) & = & \left(D_{1}w_{xx}+\lambda_{1}w(1-w)\right)I_{[M_{inact}]}(m)\nonumber \\
 &  & +\left(D_{2}w_{xx}+\lambda_{2}w(1-w)\right)I_{[M_{act}]}(m).\label{eq:fisher_two_subpop}
\end{eqnarray}
An explicit form for \eqref{eq:fisher_two_subpop} thus requires determining
how much of the population is in the active and inactive populations
over time. This is determined with the activation curves by calculating

\begin{eqnarray}
h(t;\underline{m})<m_{crit} & \iff & F(t)<\sigma(m_{crit};\underline{m})\nonumber \\
 & \iff & \underline{m}<\sigma^{-1}\left(-F(t);m_{crit}\right)=:\psi(t).\label{eq:psi}
\end{eqnarray}
Thus, $\underline{m}=\sigma^{-1}(-F(t);m)$ maps the distribution
along $m$ at time $t$ back to the initial distribution, $\phi_{1}(\underline{m}),$
and $\psi(t)$ denotes the threshold value in $\underline{m}$ between
the active and inactive populations over time. $\Phi_{1}(\psi(t))$
thus denotes the portion of the population in the inactive population,
and $1-\Phi_{1}(\psi(t))$ denotes the portion in the active population
over time.

We thus derive a nonautonomous PDE for $w$, which we will term the
\emph{averaged nonautonomous Fisher's Equation}, as: 
\begin{eqnarray}
w_{t} & = & D(t)w_{xx}+\lambda(t)w(1-w),\label{eq:w_t_phenotype}\\
w(t=0,x) & = & \phi_{2}(x)\nonumber \\
w(t,x=-\infty)=1 &  & w(t,x=\infty)=0\nonumber 
\end{eqnarray}
where 
\begin{eqnarray*}
D(t) & = & D_{2}+(D_{1}-D_{2})\Phi_{1}(\psi(t))\\
\lambda(t) & = & \lambda_{2}+(\lambda_{1}-\lambda_{2})\Phi_{1}(\psi(t)).
\end{eqnarray*}

\subsection{Three biologically-motivated examples\label{sub:Three-biologically-motivated-exa}}

We next consider three examples of \eqref{eq:nonautonomous_fisher_eqn}
that pertain to common patterns of biochemical activity during wound
healing. We will use numerical simulations of \eqref{eq:w_t_phenotype}
to investigate how different patterns of activation and deactivation
over time affect the averaged cell population profile. We will also
investigate how the profile changes when crossing the activation and
entire activation thresholds derived in \eqref{eq:transition_criterion}
and \eqref{eq:full_transition_criterion}. In each example, we fix
$m_{crit}=0.5,D_{1}=0.01,D_{2}=1,\lambda_{1}=.25$, $\lambda_{2}=0.0025,\phi_{1}(m)=\nicefrac{10}{3}I_{[.05,0.25]}(m),\phi_{2}(x)=I_{(-\infty,5]}(x)$
and $g(m)=\alpha m(1-m),$ and use a different terms for $f(t)$ to
mimic different biological situations. The choice for $g(m)$ ensures
that the distribution along $m$ stays between $m=0$ and $m=1$.
A standard central difference scheme is used for numerical simulations
of \eqref{eq:w_t_phenotype}.

\subsubsection*{Example 1: Single Sustained MAPK activation wave: $f(t)=1$\label{exa:no_forcing_threshold}}

In this example, we consider a case where we observe the entire cell
population approach a level of $m=1$ over time. Such a scenario may
represent the sustained wave of ERK 1/2 activity observed in MDCK
cells from \citep{matsubayashi_erk_2004}. The authors of \citep{posta_mathematical_2010}
proposed that the autocrine production of EGF caused this activation
in the population. We use $f(t)=1$ to observe this behavior. 

Using \eqref{eq:sigma} and \eqref{eq:psi}, we find
\begin{eqnarray*}
\sigma(m;\underline{m}) & = & \frac{1}{\alpha}\log\left(\frac{m}{1-m}\frac{1-\underline{m}}{\underline{m}}\right);\ \underline{m},m\in(0,1)\\
h(t;\underline{m}) & = & \sigma^{-1}(t;\underline{m})\ =\ \underline{m}\left((1-\underline{m})e^{-\alpha t}+\underline{m}\right)^{-1}\\
\psi(t) & = & \mbox{\ensuremath{\left(1+e^{\alpha t}\right)}}^{-1}
\end{eqnarray*}
These functions demonstrate that the distribution along $m$ is always
activating along $m$ but never reaches the $m=1$ line, as $\sigma(m;\underline{m})\rightarrow\infty$
as $m\rightarrow1^{-}$ for any $\underline{m}\in(0,1).$ The entire
population (excluding $\underline{m}=0)$ approaches $m=1$ asymptotically,
however, as $\lim_{t\rightarrow\infty}\sigma^{-1}(t;\underline{m})=1$.
In Figure \ref{fig:p(t,m)_ex1}, we use \eqref{eq:u(tm)_soln_exact_fs}
to depict the activation profile, $p(t,m)$, over time to show the
activation behavior of the population. As expected, we observe the
entire population converging to $m=1$. We include some specific plots
of the activation curves, $h(t;\underline{m}),$ for this example.
Note that the density changes along these curves by the height function
$\frac{g(\sigma^{-1}(-F(t),m))}{g(m)}$, which is equivalent to the
height function of the self-similar traveling wave ansatz made in
\eqref{eq:self_similar_TW_ansatz}. 

\begin{figure}
\centering{}\includegraphics[width=0.45\textwidth]{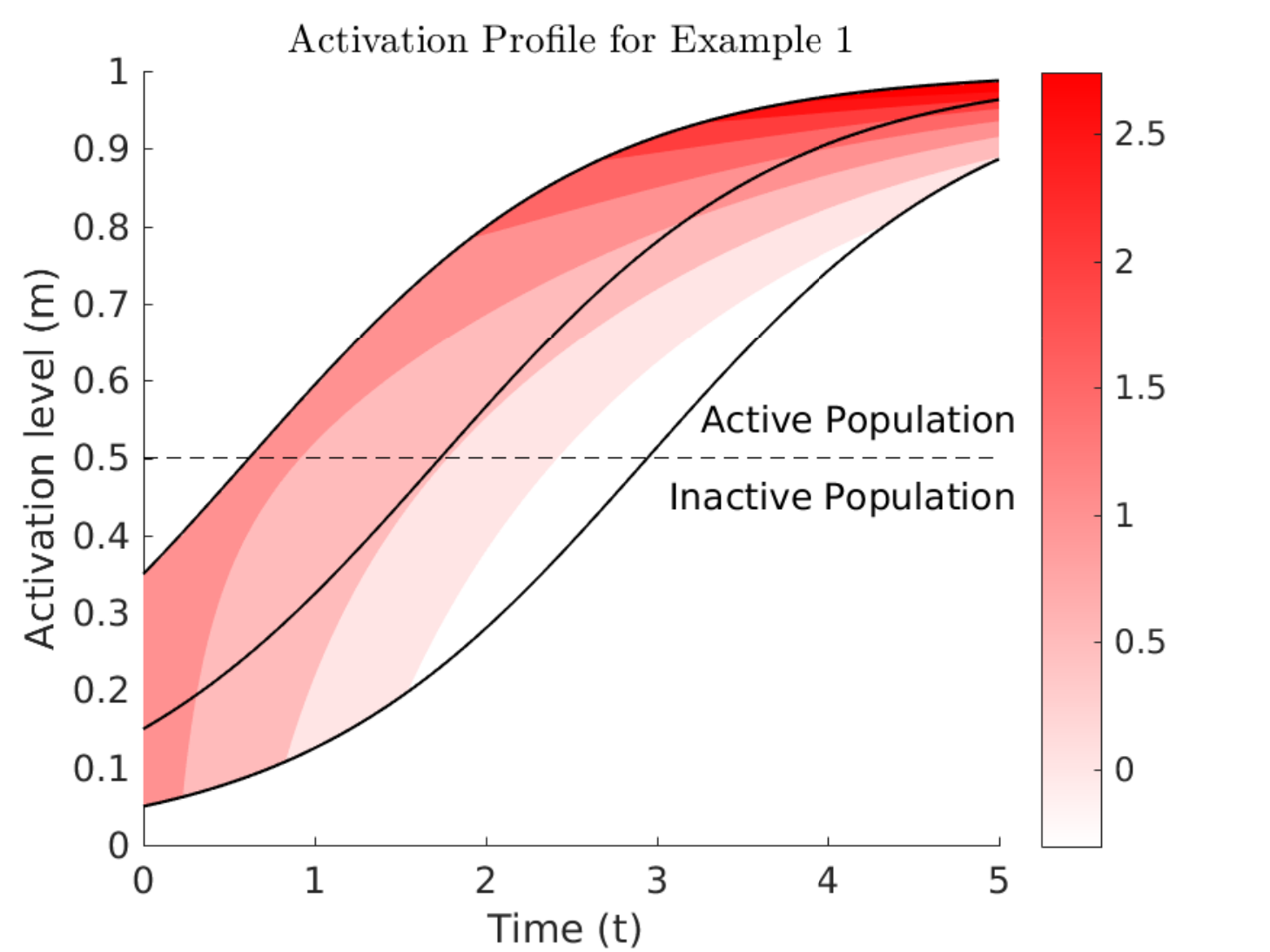}\protect\caption{The analytical solution for the activation profile, $p(t,m),$ for
Example 1 for $\alpha=0.5,$ and $\phi_{1}(m)=I_{(0.05,0.35)}(m).$
The solid black curves denote $h(t;\underline{m})$ for $\underline{m}=0.05,0.15,$
and $0.35$ and the dashed line denotes $m=m_{crit}.$ Note that a
log scale is used along $p$ for visual ease. \label{fig:p(t,m)_ex1}}
\end{figure}

\begin{figure}
\includegraphics[width=0.45\textwidth]{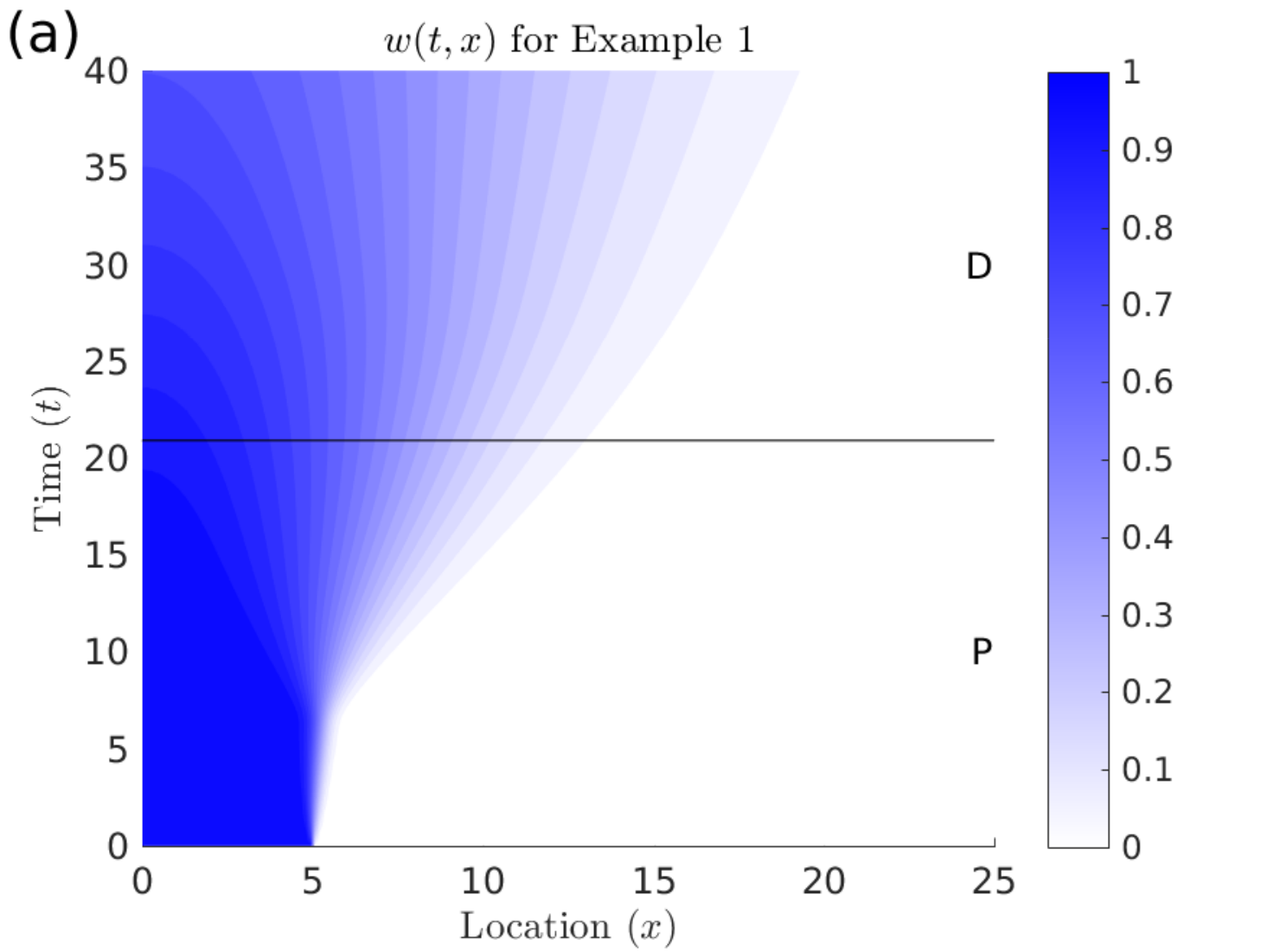}
\hfill{}\includegraphics[width=0.45\textwidth]{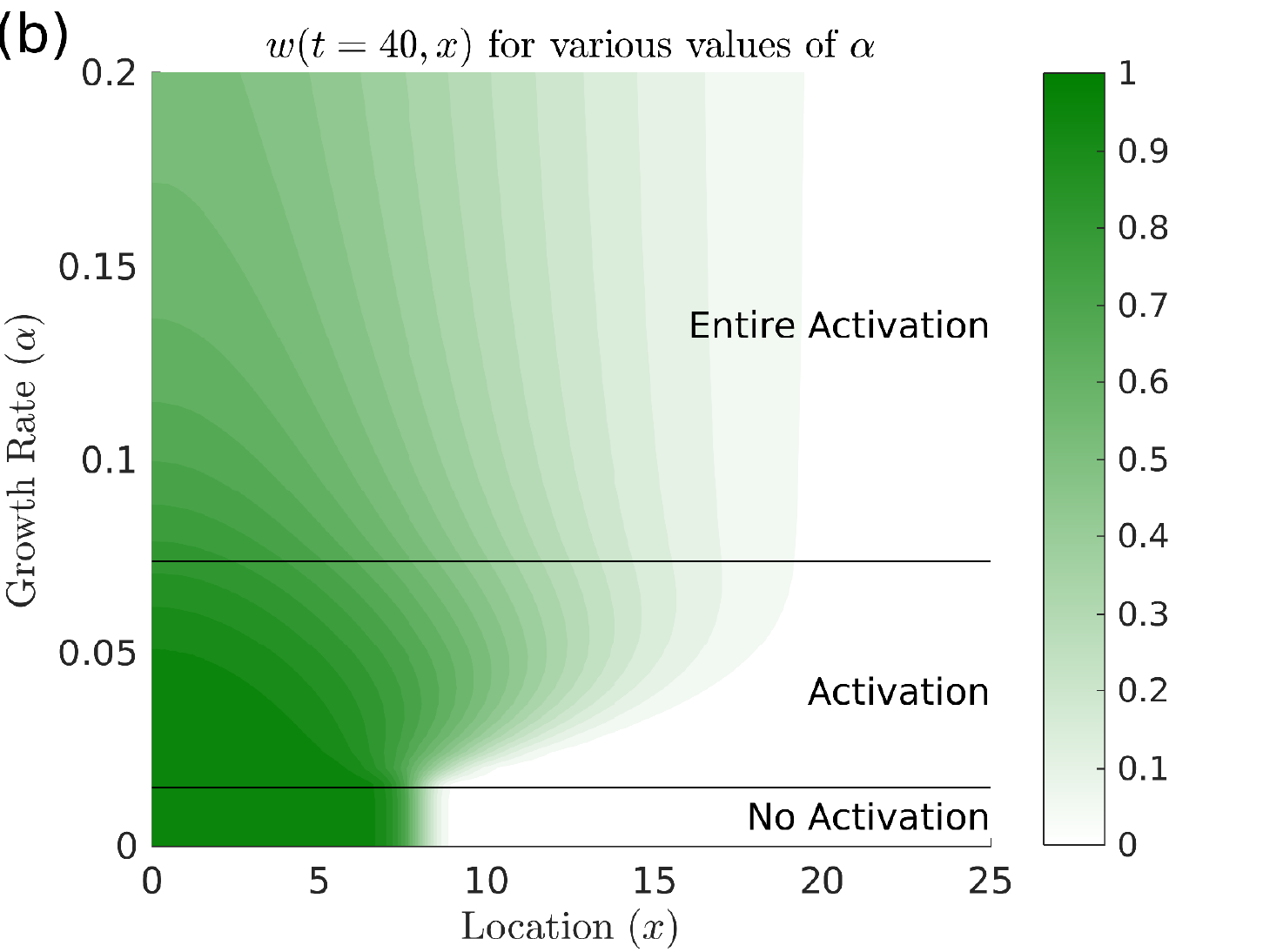}

\protect\caption{Numerical simulations of the averaged nonautonomous Fisher's equation
for Example 1. In (a), we depict a simulation of $w(t,x)$ over time
for $\alpha=0.05$. The letters ``P'' and ``D'' denote when the
population is primarily proliferating or diffusing, respectively.
In (b), we depict how the profile for $w(t=40,x)$ changes for various
values of $\alpha$. The descriptions ``No Activation'', ``Activation'',
and ``Entire Activation'' denote values of $\alpha$ for which the
population is entirely in the inactive population, split between the
active and inactive populations, or entirely in the active population
at $t=40$, respectively.\label{fig:ex1_nonaut_sims}}
\end{figure}

In Figure \ref{fig:ex1_nonaut_sims}(a), we depict a numerical simulation
of $w(t,x)$ over time using \eqref{eq:w_t_phenotype}. The slices
denoted as ``P'' and ``D'' denote when the population is primarily
proliferating ($\Phi_{1}(\psi(t))>1/2$) or diffusing $(\Phi_{1}(\psi(t))\le1/2)$
over time. The profile maintains a high cell density but limited migration
into the wound during the proliferative phase and then migrates into
the wound quickly during the diffusive phase but can not maintain
a high cell density throughout the population. In Figure \ref{fig:ex1_nonaut_sims}(b),
we investigate how the profile of $w(t=40,x)$ changes as $\alpha$
varies from $\alpha=0$ to $\alpha=0.2$. In the slice denoted ``No
activation'', the entire population is still in the inactive population
at $t=40$ and thus does not progress far into the wound or change
with $\alpha$. In the slice denoted ``Activation,'' the population
is split between the active and inactive populations at $t=40.$ The
profiles here are sensitive to increasing values of $\alpha$, as
they migrate further into the wound while maintaining a high density
near $x=0$. The slice denoted as ``Entire Activation'' denotes
simulations that are entirely in the active population by $t=40.$
As $\alpha$ increases, these simulations do not migrate much further
into the wound but do have decreasing densities at $x=0.$ These results
suggest that a combination of proliferation and diffusion must be
used to maximize population migration while maintaining a high cellular
density behind the population front. The optimal combination appears
to occur at the entire activation threshold.

\subsubsection*{Example 2: Single pulse of MAPK activation: $f(t)=\beta e^{\gamma t}-1$}

We now detail an example that exhibits a pulse of activation in the
$m$ dimension, which may represent the transient wave of ERK 1/2
activation observed in MDCK cells in \citep{matsubayashi_erk_2004}.
The authors of \citep{posta_mathematical_2010} proposed that this
wave may be caused by the rapid production of ROS in response to the
wound, followed by the quick decay of ROS or its consumption by cells.
We now let $f(t)=\beta e^{\gamma t}-1$. This forcing function arises
if ROS is present but decaying exponentially over time and modeled
by $s(t)=\beta e^{\gamma t},\beta>0,\gamma<0$ and cells activate
linearly in response to the presence of ROS but have a baseline level
of deactivation, which may be given by $f(s)=s-1.$

We see that $\sigma(m;\underline{m})$ and $\sigma^{-1}(t;\underline{m})$
are the same as in Example 1 and now compute 
\begin{eqnarray*}
h(t;\underline{m}) & = & \underline{m}\left(\underline{m}+(1-\underline{m})\exp\left[\alpha t-\frac{\alpha\beta}{\gamma}(\exp(\gamma t)-1)\right]\right)^{-1}\\
\psi(t) & = & \left(1+\exp\left[-\alpha t+\frac{\alpha\beta}{\gamma}(\exp(\gamma t)-1)\right]\right)^{-1}.
\end{eqnarray*}
In Figure \ref{fig:p(t,m)_ex2}, we use \eqref{eq:u(tm)_soln_exact_fs}
to depict the activation profile, $p(t,m),$ over time to show the
activation behavior of the population. We also include some specific
plots of the activation curves, $h(t;\underline{m})$, which show
a pulse of MAPK activity in the population that starts decreasing
around $t=5.$ Note that $h(t;0.35)$ crosses the $m=m_{crit}$ line
but $h(t;0.05)$ does not, so \eqref{eq:transition_criterion} is
satisfied for this parameter set (the population becomes activated)
but \eqref{eq:full_transition_criterion} is not (the entire population
does not become activated).

\begin{figure}
\centering{}\includegraphics[width=0.45\textwidth]{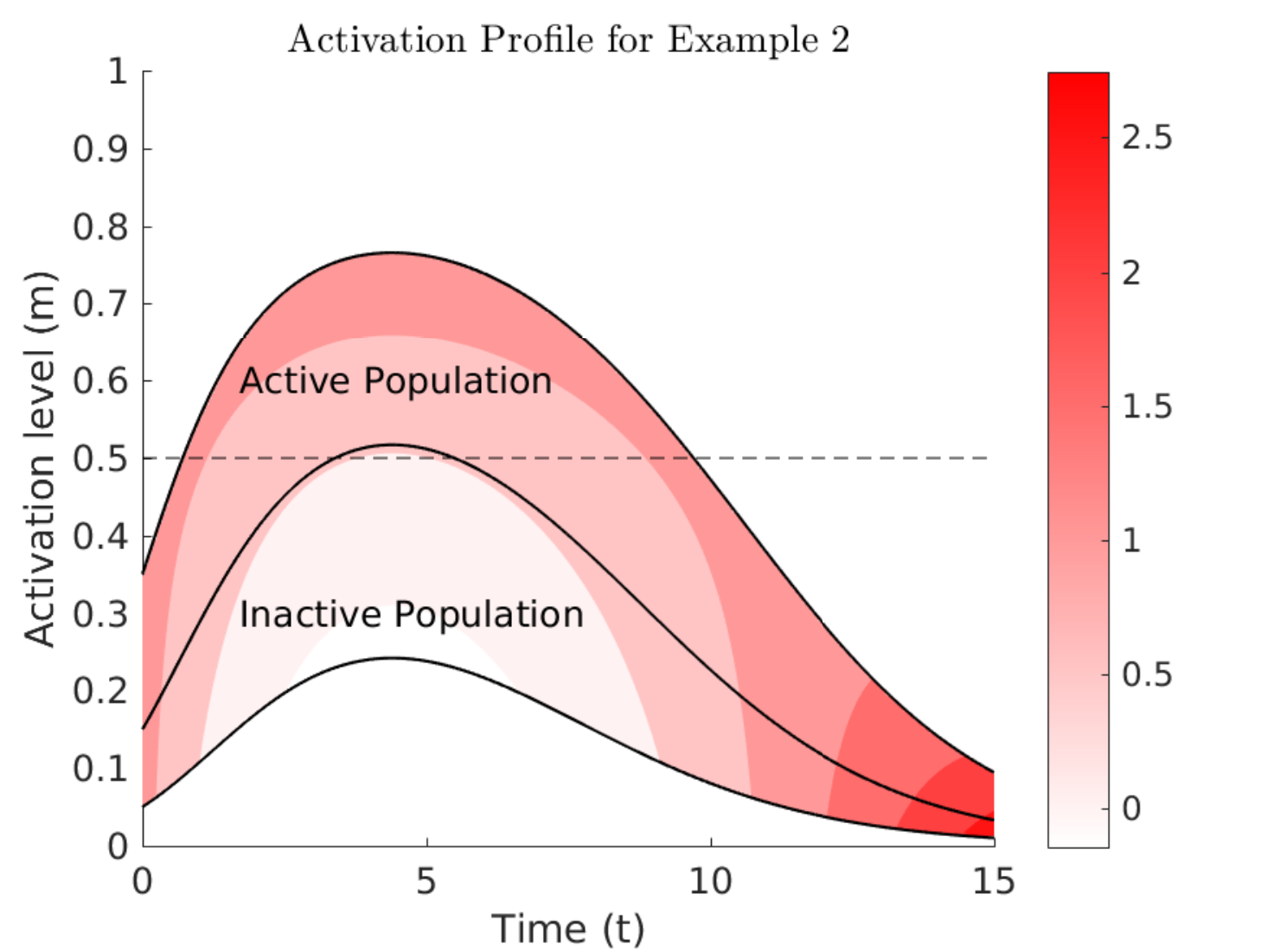}\protect\caption{The analytical solution for the activation profile, $p(t,m)$, for
Example 2 for $\alpha=0.5,\beta=3,\gamma=-1/4$ and $\phi_{1}(m)=I_{(0.05,0.35)}(m).$
The solid black curves denote $h(t;\underline{m})$ for $\underline{m}=0.05,0.15,$
and $0.35$ and the dashed line denotes $m=m_{crit}.$ Note that a
log scale is used along $p$ for visual ease. \label{fig:p(t,m)_ex2}}
\end{figure}

Using \eqref{eq:transition_criterion}, we determine our activation
criterion for this example as
\[
\frac{1-\beta+\log\beta}{\gamma}>\frac{1}{\alpha}\log\left(\frac{m_{crit}}{1-m_{crit}}\frac{1-\underline{m}_{\max}}{\underline{m}_{\max}}\right).
\]
If we fix $\gamma=-1,\alpha=1,m_{crit}=0.5,\underline{m}_{\max}=0.35,$
and $\underline{m}_{\min}=0.05$, we find that the above inequality
is satisfied for $\beta$ approximately greater than 2.55. This may
represent a scenario in which we know the decay rate of the ROS through
$\gamma,$ the activation rate of the MAPK signaling cascade through
$\alpha$, the MAPK activation distribution before ROS release with
$\underline{m}_{\min}\mbox{ and }\underline{m}_{max},$ and the activation
threshold with $m_{crit}.$ The values of $\beta$ denote the concentration
of released ROS, which should be at least 2.55 to see the population
activate. We similarly find that the entire population will activate
at some time for $\beta>5.68.$

\begin{figure}
\includegraphics[width=0.45\textwidth]{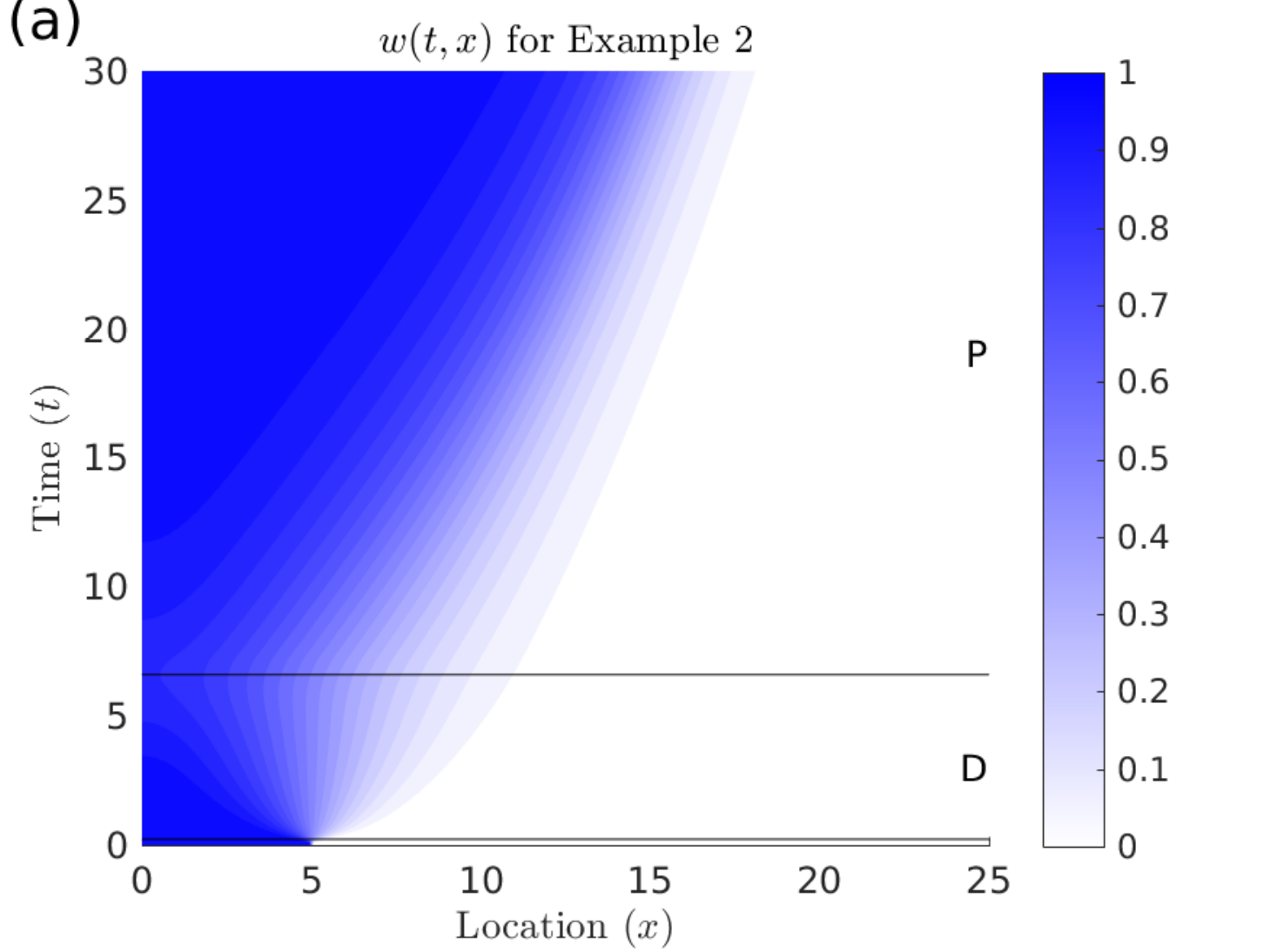}\hfill{}\includegraphics[width=0.45\textwidth]{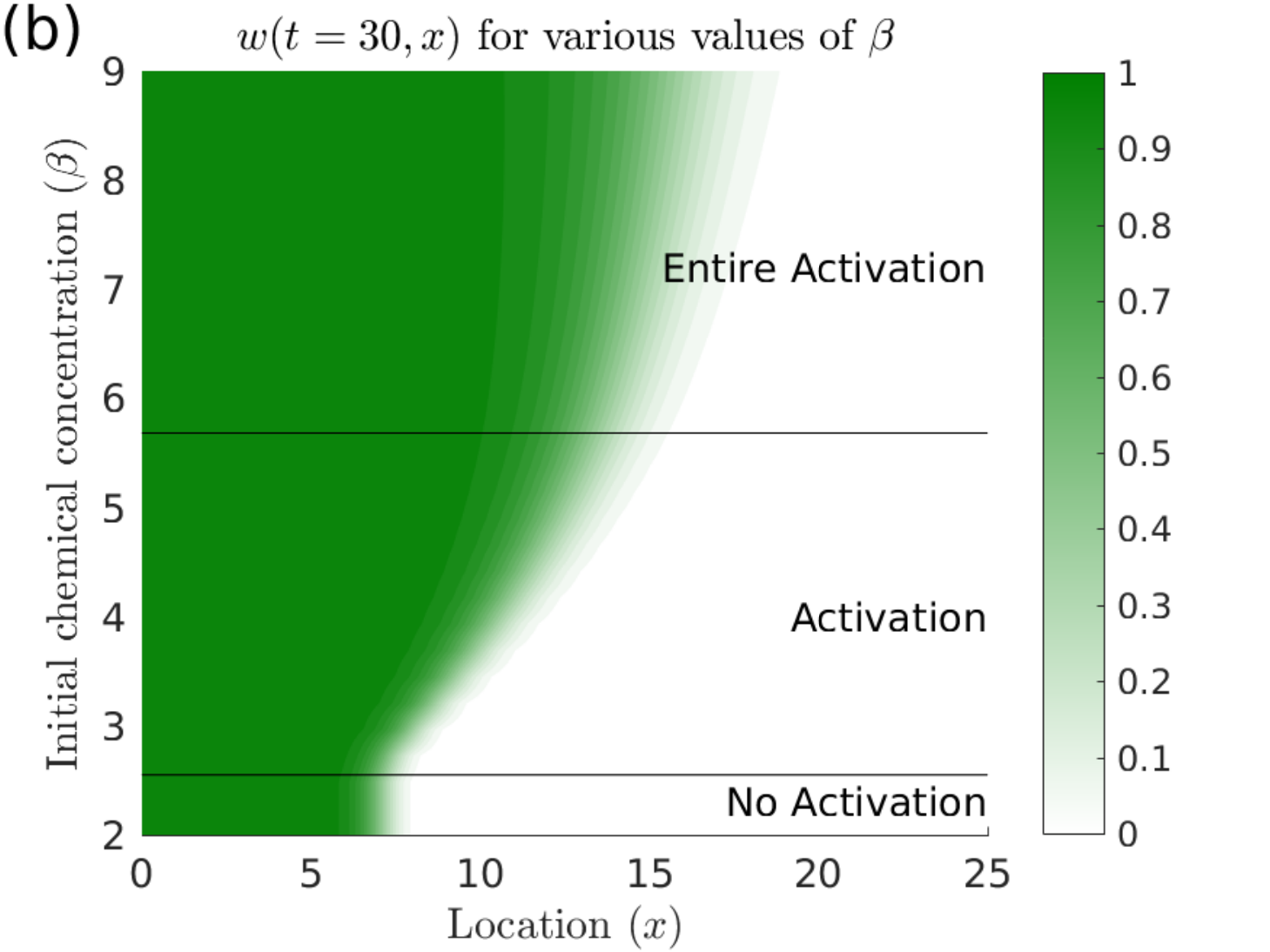}\protect\caption{Numerical simulations of the averaged nonautonomous Fisher's equation
for Example 2. In (a), we depict a simulation of $w(t,x)$ over time
for $\alpha=1,\beta=8,\gamma=-1$. Slices denoted with a ``P'' or
``D'' denote when the population is primarily proliferating or diffusing,
respectively. In (b), we depict how the profile for $w(t=30,x)$ changes
for various values of $\beta$. The descriptions ``No activation'',
``Activation'', and ``Entire Activation'' denote values of $\beta$
for which the population is entirely in the inactive population, split
between the active and inactive populations, or entirely in the active
population at $t=t_{\max}$. \label{fig:ex2_nonaut_sims}}
\end{figure}

In Figure \ref{fig:ex2_nonaut_sims}(a), we depict a numerical simulation
of \eqref{eq:w_t_phenotype} for this example. The population quickly
transitions to a diffusing stage due to the pulse of MAPK activation
and shows the smaller densities ($u$ approximately less than $0.2)$
migrating into the wound rapidly while the density behind the population
front drops. As the pulse of MAPK activation ends and the population
transitions back to a proliferating phenotype, the populations restores
a high density behind the cell front and begins to develop a traveling
wave profile, as suggested by the parallel contour lines. In Figure
\ref{fig:ex2_nonaut_sims}(b), we investigate how the profile for
$w(t=30,x)$ changes as $\beta$ varies from $\beta=2$ to $\beta=9$
while keeping all other parameters fixed. We observe that the profile
is the same for all values of $\beta<2.55,$ as \eqref{eq:transition_criterion}
is not satisfied. As $\beta$ increases past the activation threshold,
the profile shows increased rates of migration into the wound. After
passing the entire activation threshold \eqref{eq:full_transition_criterion},
the profile continues to migrate further as $\beta$ increases, but
appears less sensitive to $\beta$. This increased migration is likely
due to the population spending more time in the active population
for larger values of $\beta$. Note that for all simulations shown,
the pulse of MAPK activation has finished by $t=30.$

\subsubsection*{Example 3: Periodic pulses of MAPK activation: $f(t)=\beta\sin(\gamma t)$}

As a last example, we exhibit a scenario with periodic waves of activity.
Such behavior was observed in some of the experiments performed in
\citep{zi_quantitative_2011}, in which cell cultures of the HaCaT
cell line were periodically treated with TGF-$\beta$ to investigate
how periodic treatment with TGF-$\beta$ affects activation of the
SMAD pathway (the canonical pathway for TGF-$\beta$, which also influences
cell proliferation and migration). We let $f(t)=\beta\sin(\gamma t),\beta,\gamma>0$,
which occurs if the concentration of TGF-$\beta$ over time is given
by $s(t)=1+\sin(\gamma t),$ and cells activate linearly in response
to $s$ and have a baseline rate of deactivation, given by $f(s)=\beta(s-1)$. 

We now calculate

\[
h(t;\underline{m})=\underline{m}\left(\underline{m}+(1-\underline{m})\exp\left[\frac{\alpha\beta}{\gamma}\left(\cos(\gamma t)-1\right)\right]\right)^{-1}
\]

\[
\psi(t)=\left(1+\exp\left[\frac{\alpha\beta}{\gamma}(1-\cos(\gamma t))\right]\right)^{-1}.
\]
In Figure \ref{fig:p(t,m)_ex3}, we use \eqref{eq:u(tm)_soln_exact_fs}
to depict the activation profile, $p(t,m),$ over time to show the
activation behavior of the population. We also include some specific
plots of the activation curves $h(t;\underline{m})$, which demonstrate
periodic waves of activation along $m$. Note that $h(t;0.05)$ crosses
the $m=m_{crit}$ line, so \eqref{eq:full_transition_criterion} is
satisfied, and the entire population becomes activated at some points
during the simulation.

\begin{figure}
\centering{}\includegraphics[width=0.45\textwidth]{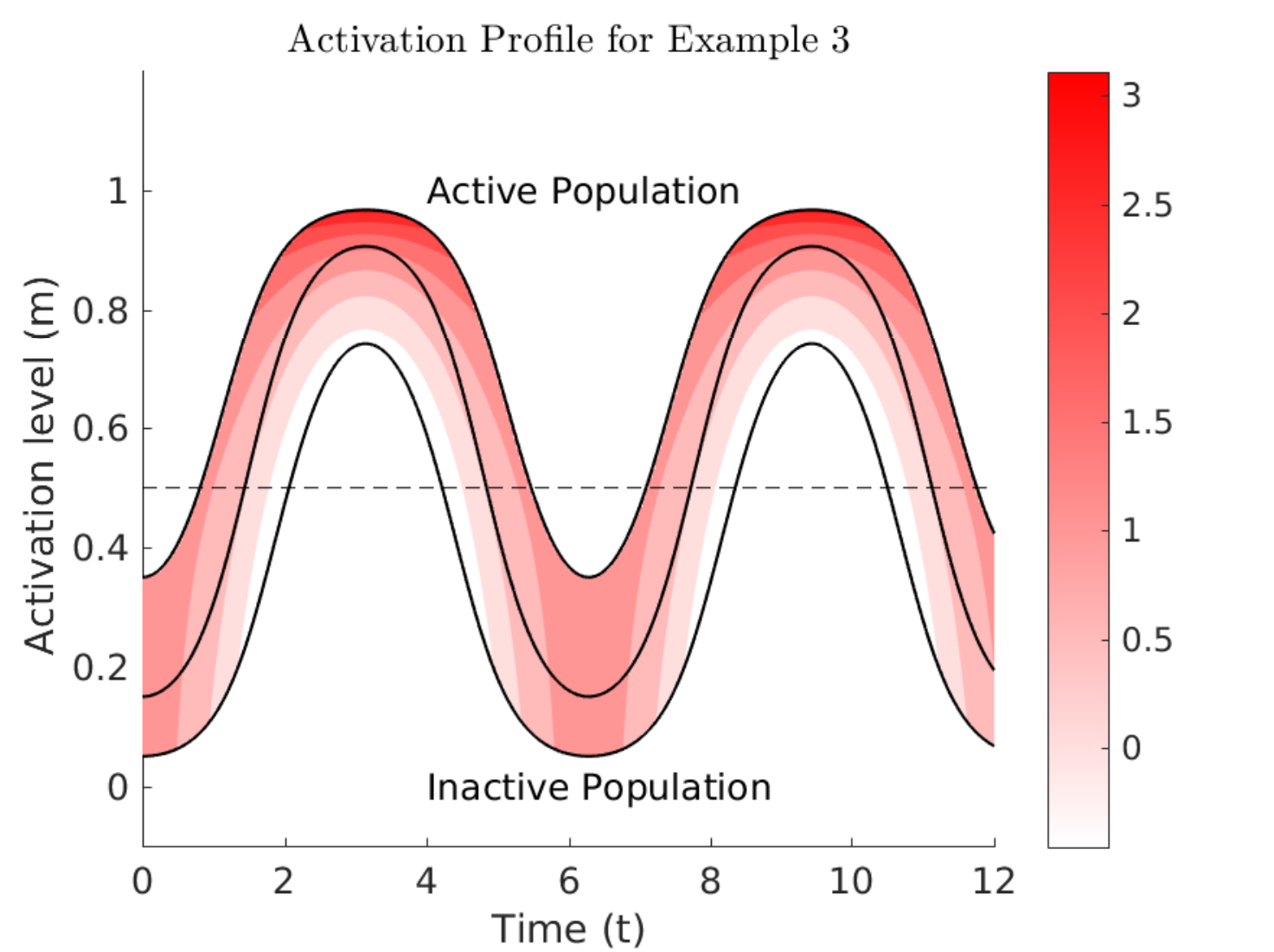}\protect\caption{The analytical solution for the activation profile, $p(t,m)$, for
Example 3 for $\alpha=1/2,\beta=4,\gamma=1$ and $\phi_{1}(m)=I_{(0.05,0.35)}(m).$
The solid black curves denote $h(t;\underline{m})$ for $\underline{m}=0.05,0.15,$
and $0.35$ and the dashed line denotes $m=m_{crit}.$ Note that a
log scale is used along $p$ for visual ease. \label{fig:p(t,m)_ex3}}
\end{figure}

The activation criterion \eqref{eq:transition_criterion} can be solved
as
\[
\frac{2\beta}{\gamma}>\frac{1}{\alpha}\log\left(\frac{m_{crit}}{1-m_{crit}}\frac{1-\underline{m}_{\max}}{\underline{m}_{\max}}\right).
\]
We thus calculate that if we fix $\beta=1,\alpha=1/2,\underline{m}_{max}=0.35,\underline{m}_{\min}=0.05,$
and $m_{crit}=0.5,$ then the activation criterion \eqref{eq:transition_criterion}
is satisfied for $\gamma<1.615$ and the entire activation criterion
\eqref{eq:full_transition_criterion} is satisfied for $\gamma<0.34.$
These estimates would tell us how frequently signaling factor treatment
is needed to see different patterns of activation in the population.

\begin{figure}
\includegraphics[width=0.45\textwidth]{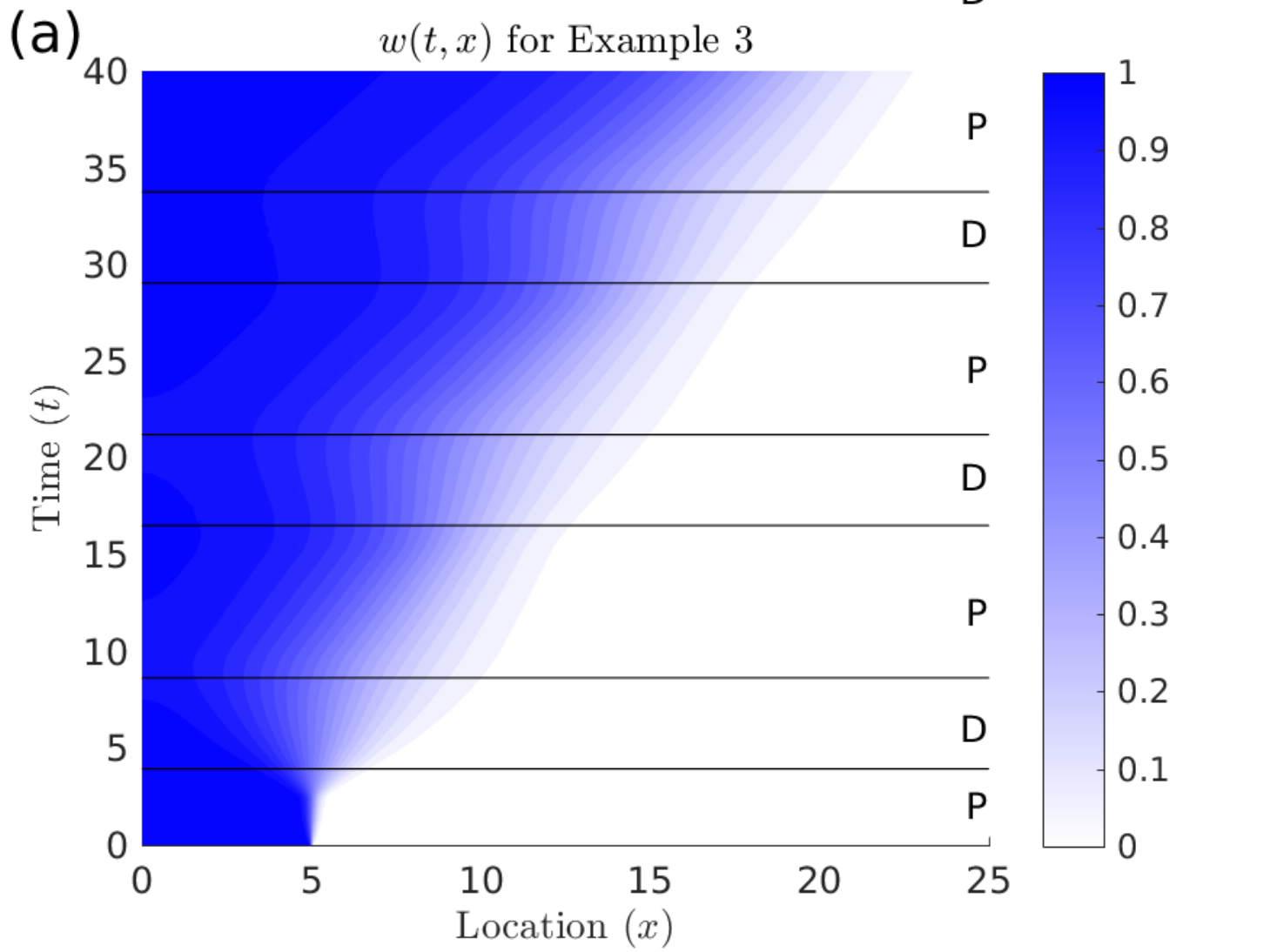}\hfill{}\includegraphics[width=0.45\textwidth]{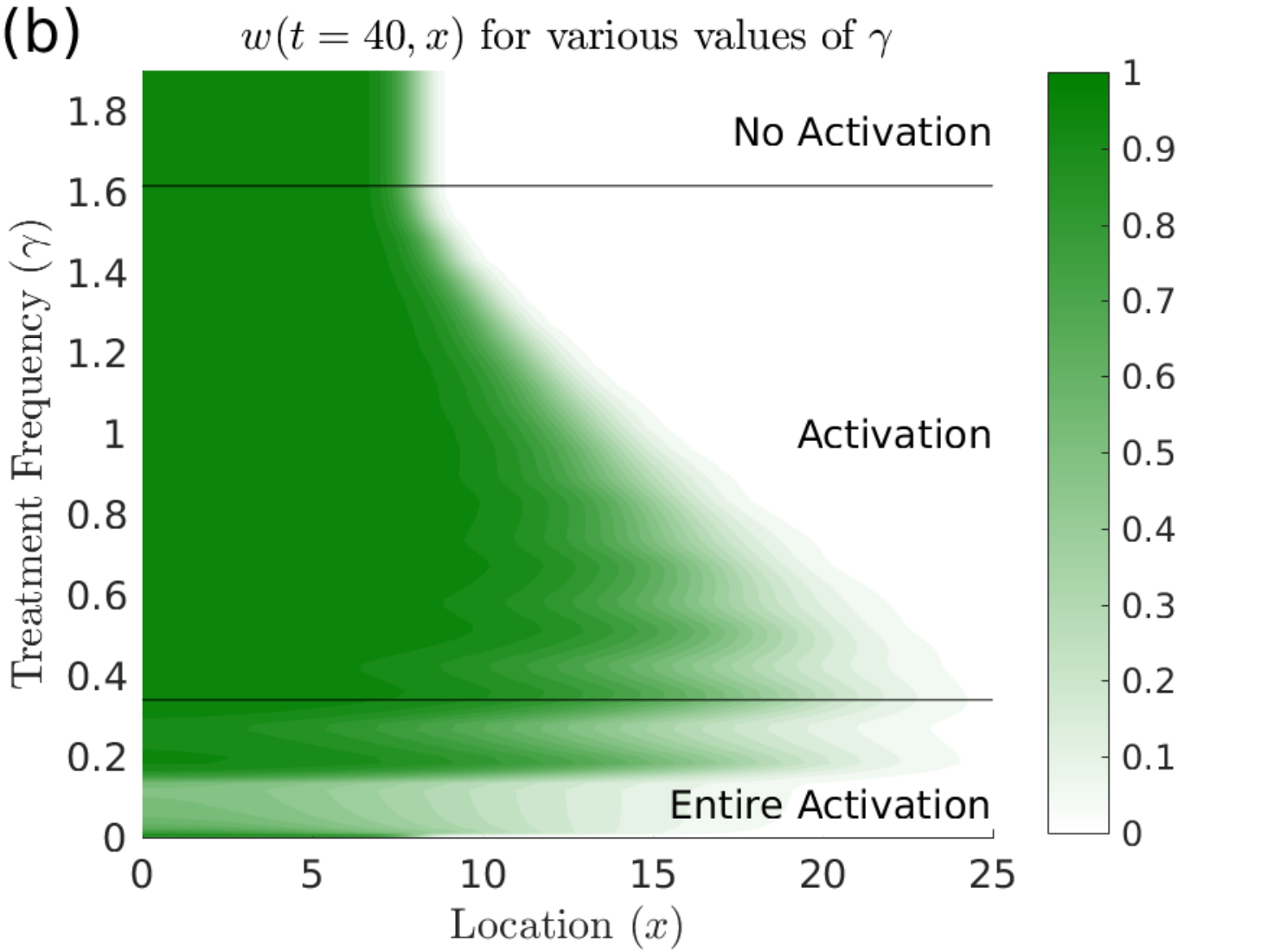}\protect\caption{Numerical simulations of the averaged nonautonomous Fisher's equation
for Example 3. In (a), we depict a simulation of $w(t,x)$ over time
for $\alpha=0.5,\beta=1,$ and \protect \\
$\gamma=1/2$. Slices denoted with a ``P'' or ``D'' denote when
the population is primarily proliferating or diffusing, respectively.
In (b), we depict $w(t=40,x)$ for various values of $\gamma$. The
descriptions ``No activation'', ``Activation'', and ``Entire
Activation'' denote values of $\gamma$ for which the population
is entirely in the inactive population, split between the active and
inactive populations, or entirely in the active population at $t=t_{\max}$.\label{fig:ex3_nonaut_sims}}
\end{figure}

In Figure \ref{fig:ex3_nonaut_sims}(a), we depict a numerical simulation
of \eqref{eq:w_t_phenotype} for this example. The population phenotype
has a period of $4\pi$, and we see that the lower densities migrate
into the wound most during the diffusive stages, whereas all densities
appear to migrate into the wound at similar speeds during the proliferative
stages. In Figure \ref{fig:ex3_nonaut_sims}(b), we investigate how
the profile for $w(t=40,x)$ changes as $\gamma$ varies between $\gamma=0$
and $\gamma=1.9$ while keeping all other parameters fixed. All profiles
appear the same for $\gamma>1.615$ as \eqref{eq:transition_criterion}
is not satisfied. As $\gamma$ decreases below this threshold, more
of the population becomes activated during the simulation, culminating
in a maximum propagation of the population at the entire activation
threshold, $\gamma\approx0.34.$ As $\gamma$ falls below $\mbox{\ensuremath{\gamma}=}0.34,$
the population tends to migrate less, although the population does
migrate far for $\gamma$ near $0.2$. For $\gamma<0.2,$ the population
appears to spend too much time in the active population and diffuses
excessively with limited proliferation. These simulations lead to
shallow profiles that do not migrate far into the wound. As $\gamma$
approaches zero, the simulations would become entirely activated,
but do not before $t=40$. These simulations stay in the inactive
population for the duration of the simulation and do not migrate far
into the wound. 

\section{Discussion and Future work\label{sec:Discussion-and-future}}

We have investigated a structured Fisher's Equation that incorporates
an added dimension for biochemical activity that influences population
migration and proliferation. The method of characteristics proved
to be a useful way to track the progression along the population activity
dimension over time. With the aid of a phase plane analysis and an
asymptotically autonomous Poincare-Bendixson Theorem, we were able
to prove the existence of a self-similar traveling wave solution to
the equation when diffusion and proliferation do not depend on MAPK
activity. The height function of the self-similar traveling wave ansatz
along characteristic curves is demonstrated in Figures \ref{fig:p(t,m)_ex1},
\ref{fig:p(t,m)_ex2}, and \ref{fig:p(t,m)_ex3}. We believe our analysis
could be extended to investigate structured versions of other nonlinear
PDEs.

Activation of the MAPK signaling cascade is known to influence collective
migration during woung healing through cellular migration and proliferation
properties. For this reason, we also considered a structured PDE model
in which the rates of cellular diffusion and proliferation depend
on the levels of MAPK activation in the population. We also extended
the model to allow for the presence of an external cytokine or growth
factor that regulates activation and deactivation along the MAPK signaling
cascade. We derived two activation criteria for the model to establish
conditions under which the population will become activated during
simulations. As numerical simulations of the structured equation are
prone to error via numerical diffusion, we derived a nonautonomous
equation in time and space to represent the average population behavior
along the biochemical activity dimension. Using this nonautonomous
equation, we exhibited three simple examples that demonstrate biologically
relevant activation levels and their effects on population migration:
a sustained wave of activity, a pulse of activity, and periodic pulses
of activity. We found that the population tends to migrate farthest
while maintaining a high cell density at the entire activation threshold
value, \eqref{eq:full_transition_criterion}, for the sustained wave
and periodic pulse patterns of activation. The single pulse case continued
migrating further into the wound after passing the entire activation
threshold but appeared less sensitive after doing so.

A natural next step for this analysis is to use a structured population
model of this sort in combination with biological data to thoroughly
investigate the effects of MAPK activation and deactivation on cell
migration and proliferation during wound healing. Previous mathematical
models have focused on either collective migration during wound healing
assays in response to EGF treatment (while neglecting the MAPK signaling
cascade) \citep{johnston_estimating_2015,nardini_modeling_2016} or
MAPK propagation during wound healing assays (while neglecting cell
migration) \citep{posta_mathematical_2010}. To the best of our knowledge,
no mathematical models have been able to reliably couple signal propagation
and its effect on cell migration during wound healing. The examples
detailed in this work intentionally used the simplest terms possible
as a means to focus on the underlying mathematical aspects. With a
separate in-depth study into the biochemistry underlying the MAPK
signaling cascade and its relation with various cytokines or growth
factors, more complicated and biologically relevant terms for $g(m),f(s),$
and $s(t)$ can be determined to help elucidate the effects of MAPK
activation on cell migration during wound healing.

The analytical techniques used in this study cannot be used to investigate
spatial patterns of biochemical activity due to the parabolic nature
of \eqref{eq:structured_fisher_eqn} in space. Cell populations also
migrate via chemotaxis during wound healing, in which cells migrate
up a concentration gradient of some chemical stimulus \citep{ai_reaction_2015,keller_traveling_1971,landman_diffusive_2005,newgreen_chemotactic_2003}.
Chemotactic equations are hyperbolic in space, which may facilitate
spatial patterns of MAPK activation during wound healing, such as
those described experimentally in \citep{chapnick_leader_2014}. As
various pathways become activated and cross-talk during wound healing
to influence migration \citep{guo_signaling_2009}, future studies
could also investigate a population structured along multiple signaling
pathways, $u(t,x,\vec{m})$ for the vector $\vec{m}=(m_{1},m_{2},\dots,m_{n})^{T}$.
Because the cell population also produces cytokines and growth factors
for paracrine and autocrine signaling during wound healing, these
models would also benefit from unknown variables representing ROS,
TGF-$\beta$, EGF, etc.

While the main motivation for this study is epidermal wound healing,
there are potential applications in other areas of biology. Fisher's
equation has also been used to study population dynamics in ecology
and epidemiology \citep{ai_travelling_2005,hastings_spatial_2005,shigesada_biological_1997}.
Our framework could be extended to a case where an environmental effect,
such as seasonal forcing, impacts species migration or susceptibility
of individuals to disease. The results presented here may thus aid
in a plethora of mathematical biology studies.

\bibliographystyle{my_siam}
\bibliography{15_home_john_fall_2016_stucture_manuscript_SIAM_final_MyLibrary}

\appendix

\section{Properties of $\sigma^{-1}(t;\underline{y})$\label{sec:Properties-of}}

If we assume that $g$ is positive and uniformly continuous, then
$\sigma^{-1}(t;\underline{y})$ exists and satisfies the following:
\begin{equation}
\frac{d}{dt}\sigma^{-1}(t;\underline{y})=g(\sigma^{-1}(t;\underline{y})),\ \sigma^{-1}(0;\underline{y})=\underline{y}.\label{eq:dsigma_deriv}
\end{equation}
 To derive \eqref{eq:dsigma_deriv}, see that
\begin{eqnarray*}
y(t) & = & \sigma^{-1}(t;\underline{y})\\
\Rightarrow\sigma(y(t);\underline{y}) & = & t\\
\Rightarrow\frac{d}{dt}\left(\sigma(y(t);\underline{y})\right)=\frac{d}{dy}\sigma(y(t);s)\frac{dy}{dt} & = & 1\\
\Rightarrow\frac{1}{g(y(t))}\frac{dy}{dt} & = & 1\\
\Rightarrow\frac{dy}{dt} & = & g(y(t))\\
\Rightarrow\frac{d}{dt}\sigma^{-1}(t;\underline{y}) & = & g(\sigma^{-1}(t;y)).
\end{eqnarray*}
and for the initial condition, 
\begin{eqnarray*}
\sigma(\underline{y},\underline{y}) & = & 0\\
\Rightarrow\sigma^{-1}(0,\underline{y}) & = & \underline{y}.
\end{eqnarray*}

\section{Derivation of \eqref{eq:v_t_diffeq}\label{sec:Derivation-of}}

In \eqref{eq:size_characteristic}, we defined
\[
v(t;\underline{y}):=u(t,y=\sigma^{-1}(t;\underline{y})).
\]
Taking the derivative of $v(t;\underline{y})$ with respect to time,
we find with the aid of the chain rule:
\begin{eqnarray*}
\frac{d}{dt}v(t;\underline{y}) & = & \frac{\partial}{\partial t}u(t,y=\sigma^{-1}(t;\underline{y}))+\frac{\partial}{\partial y}u(t,y=\sigma^{-1}(t;\underline{y}))\cdot\frac{d}{dt}\sigma^{-1}(t;\underline{y})\\
 & = & -\frac{\partial}{\partial y}\left[g\left(\sigma^{-1}(t;\underline{y})\right)u\left(t,y=\sigma^{-1}(t;\underline{y})\right)\right]\\
 &  & +\ Au\left(t,y=\sigma^{-1}(t;\underline{y})\right)+g\left(\sigma^{-1}(t;\underline{y})\right)\frac{\partial}{\partial y}u\left(t,y=\sigma^{-1}(t;\underline{y})\right)\\
 & = & -g'(\sigma^{-1}(t;\underline{y}))v(t;\underline{y})+Av(t;\underline{y}).
\end{eqnarray*}

\section{Relevant Material on Asymptotically Autonomous Differential Systems\label{sec:Statement-of-AA_markus}}

The theorem statements in this section have been slightly modified
to match notation from our study. Consider the two vector fields,
\begin{eqnarray}
\dot{x} & = & f(t,x)\label{eq:f_nonaut}\\
\dot{y} & = & g(y),\label{eq:g_aut}
\end{eqnarray}
\ for $x,y\in\mathbb{R}^{n}$ and $t>0$. Assume $f(t,x)$ and $g(x)$
are continuous in $t$ and $x$ and locally Lipschitz in $x$ for
$x\in\Omega$ and $t>0$, where $\Omega$ is an open subset of $\mathbb{R}^{n}.$
We say that \eqref{eq:f_nonaut} is asymptotically autonomous with
limit equation \eqref{eq:g_aut} if $f(t,x)\rightarrow g(x)$ pointwise
as $t\rightarrow\infty$ on any compact subset of $\Omega$.

We will denote the $\omega$-limit sets for all points starting in
the set $\Theta\subset\Omega$ at $t=0$ for the system \eqref{eq:f_nonaut}
as $\omega_{f}(\Theta)$. This asymptotically autonomous Poincare-Bendixson
Theorem was introduced in \citep{markus_asymptotically_1956} and
states

\textbf{Asymptotically Autonomous Poincare-Bendixson Theorem}: 

Let $n=2$ and \eqref{eq:f_nonaut} be asymptotically autonomous with
limit equation \eqref{eq:g_aut} in $\Omega\subset\mathbb{R}^{2}.$
Let a solution, $x(t)$, of \eqref{eq:f_nonaut} lie in a compact
set $\Theta\subset\Omega$ for large $t$ and suppose $\omega_{f}(\Theta)$
contains no equilibria of \eqref{eq:g_aut}. Then $\omega_{f}(\Theta)$
is the union of periodic orbits of \eqref{eq:g_aut}.

The proof of this Theorem is a result of the standard Poincare-Bendixson
Theorem (see \citep[Section 6.6]{meiss_differential_2007}) and the
following theorem, which is also proved in \citep{markus_asymptotically_1956}.

\textbf{Theorem: }Let \eqref{eq:f_nonaut} be asymptotically autonomous
with limit equation \eqref{eq:g_aut} in $\Omega\in\mathbb{R}^{n}.$
Let P be a stable equilibrium point of \eqref{eq:g_aut}. Then there
is a neighborhood, $N$, of $P$ and a time $T$ such that $\omega_{f}(N)$
=$\{P\}$ for all solutions of \eqref{eq:f_nonaut} starting at time
$T$ or later.
\end{document}